\documentclass[11pt,letterpaper,english,aip,letter,superscriptaddress,nofootinbib,floatfix,jcp]{revtex4}
\usepackage{ae,aecompl}
\usepackage[T1]{fontenc}
\usepackage[latin1]{inputenc}
\usepackage{xcolor}
\usepackage{pdfcolmk}
\usepackage{amsmath}
\usepackage{amssymb}
\usepackage{graphicx}
\usepackage{esint}
\PassOptionsToPackage{normalem}{ulem}
\usepackage{ulem}

\makeatletter

\pdfpageheight\paperheight
\pdfpagewidth\paperwidth

\providecolor{lyxadded}{rgb}{1,0,0}
\providecolor{lyxdeleted}{rgb}{0,0,1}

\@ifundefined{textcolor}{}
{%
 \definecolor{BLACK}{gray}{0}
 \definecolor{WHITE}{gray}{1}
 \definecolor{RED}{rgb}{1,0,0}
 \definecolor{GREEN}{rgb}{0,1,0}
 \definecolor{BLUE}{rgb}{0,0,1}
 \definecolor{CYAN}{cmyk}{1,0,0,0}
 \definecolor{MAGENTA}{cmyk}{0,1,0,0}
 \definecolor{YELLOW}{cmyk}{0,0,1,0}
}

\usepackage{aecompl}

\renewcommand{\citet}[1]{\cite{#1}}

\usepackage{babel}

\makeatother

\usepackage{babel}
\begin{document}
\setlength{\abovedisplayskip}{0.4ex}\setlength{\belowdisplayskip}{0.4ex}
\setlength{\abovedisplayshortskip}{0.15ex}\setlength{\belowdisplayshortskip}{0.15ex}

\title{A Minimally-Resolved Immersed Boundary Model for Reaction-Diffusion
Problems}

\author{Amneet Pal Singh Bhalla}

\affiliation{Department of Mechanical Engineering, Northwestern University, Evanston,
IL 60208}

\affiliation{Courant Institute of Mathematical Sciences, New York University,
New York, NY 10012}

\author{Boyce E. Griffith}

\affiliation{Leon H. Charney Division of Cardiology, Department of Medicine, New
York University School of Medicine, New York, NY 10016}

\affiliation{Courant Institute of Mathematical Sciences, New York University,
New York, NY 10012}

\author{Neelesh A. Patankar}

\affiliation{Department of Mechanical Engineering, Northwestern University, Evanston,
IL 60208}

\author{Aleksandar Donev}

\email{donev@courant.nyu.edu}

\selectlanguage{english}%

\affiliation{Courant Institute of Mathematical Sciences, New York University,
New York, NY 10012}
\begin{abstract}
We develop an immersed-boundary approach to modeling reaction-diffusion
processes in dispersions of reactive spherical particles, from the
diffusion-limited to the reaction-limited setting. We represent each
reactive particle with a minimally-resolved ``blob'' using many
fewer degrees of freedom per particle than standard discretization
approaches. More complicated or more highly resolved particle shapes
can be built out of a collection of reactive blobs. We demonstrate
numerically that the blob model can provide an accurate representation
at low to moderate packing densities of the reactive particles, at
a cost not much larger than solving a Poisson equation in the same
domain. Unlike multipole expansion methods, our method does not require
analytically-computed Green's functions, but rather, computes regularized
discrete Green's functions on the fly by using a standard grid-based
discretization of the Poisson equation. This allows for great flexibility
in implementing different boundary conditions, coupling to fluid flow
or thermal transport, and the inclusion of other effects such as temporal
evolution and even nonlinearities. We develop multigrid-based preconditioners
for solving the linear systems that arise when using implicit temporal
discretizations or studying steady states. In the diffusion-limited
case the resulting linear system is a saddle-point problem, the efficient
solution of which remains a challenge for suspensions of many particles.
We validate our method by comparing to published results on reaction-diffusion
in ordered and disordered suspensions of reactive spheres.
\end{abstract}
\maketitle
\global\long\def\V#1{\boldsymbol{#1}}
 \global\long\def\M#1{\boldsymbol{#1}}
 \global\long\def\Set#1{\mathbb{#1}}

\global\long\def\D#1{\Delta#1}
 \global\long\def\d#1{\delta#1}

\global\long\def\norm#1{\left\Vert #1\right\Vert }
 \global\long\def\abs#1{\left|#1\right|}

\global\long\def\P{{P}}

\global\long\def\grad{\boldsymbol{\nabla}}
 \global\long\def\av#1{\langle#1\rangle}
 \global\long\def\peroxide{H_{2}O_{2}}
 \global\long\def\cU{\mathcal{U}}
 \global\long\def\cL{\mathcal{L}}
 \global\long\def\cLr{\mathcal{L}_{r}}
 \global\long\def\Cav{\av C}
 \global\long\def\Lap{\nabla^{2}}

\section{Introduction}

Over the past decade several types of synthetic chemically-propelled
micro- and nano-swimmers have been manufactured and studied in the
lab \cite{Nanomotors_Kapral,Hematites_Science}. Modeling the behavior
of suspensions of such particles requires including both the hydrodynamics
of the fluid flow \emph{and} the reaction-diffusion transport of the
chemically-reacting species around the particles. Reaction-diffusion
processes are ubiquitous in chemical engineering applications and
have been the subject of intense study on their own right, particularly
in the diffusion-limited regime in which the rate-limiting process
is the diffusion of the reactants to the reactive sites \cite{ReactionDiffusion_Torquato,ReactionDiffusion_Brady,ReactionDiffusion_Lu,FPKMC2}.
The fluid dynamics of active suspensions is a very active field itself
\cite{ActiveSuspensions}, and methods similar to the one developed
here have recently been proposed \cite{MicroSwimmers_Lushi}. Fewer
theoretical and computational studies have been done that couple the
Navier-Stokes equation for the fluid velocity with a reaction-diffusion
equation for the concentration of the reactants \cite{ReactionDiffusion_Fauci}.

Herein, we develop a model and numerical method for such reaction-diffusion
problems using an approach based on the immersed boundary method \cite{IBM_PeskinReview}.
We build on recently-developed computational methods for the hydrodynamics
of dilute suspensions of rigid spherical particles \cite{DirectForcing_Balboa,ISIBM}.
In the \emph{minimally-resolved} approach that we adopt, each particle
is represented by a single point-like object that we will call a \emph{blob}.
The blob does not resolve the details of the particle surface, but
rather, represents the effective contribution of the particle using
a small number of degrees of freedom, as do particles in (stochastic)
dissipative particle dynamics \cite{SDPD_Scaling}. This is very similar
to the use of multipole expansions to solve fluid-particle and reaction-diffusion
problems \cite{ReactionDiffusion_Brady,BrownianDynamics_OrderNlogN}.
The essential difference is that in our approach we do not employ
analytical Green's functions to represent the response of the continuum
model to the particles. Instead, we use a standard grid-based solver
for the continuum equations to generate the required response function
``on the fly''. This approach allows us to easily change the boundary
conditions, including cases in which the concentration of the chemical
reactants affects the fluid flow via osmo-phoretic effects \cite{Hematites_PRL}.
Furthermore, temporal evolution (unsteadiness) and even nonlinearities
can easily be added, and the resolution can be refined locally using
adaptive mesh refinement techniques \cite{IBAMR}.

The method developed here can be used as the core component of a more
traditional immersed-boundary method \cite{IBM_PeskinReview} for
reaction-diffusion problems, in which a large number of blobs (traditionally
called markers or regularized delta functions) are placed on the surface
at which the reactive boundary condition is specified. In fact, by
changing the number of blobs used to discretize, for example, a spherical
particle, one can go from high-resolution models (e.g., hundreds of
blobs per particle), to low-resolution models (e.g., 12 blobs per
particle), to minimally-resolved models (a single blob per particle).
In this work we focus on minimally-resolved modeling of reaction-diffusion
problems in particle suspensions, extending the traditional scope
and the traditional view of convergence of the immersed boundary method
\cite{DirectForcing_Balboa,ISIBM}. In our approach, even though each
particle is represented only by the position of its centroid, it is
not appropriate to consider it a ``point'' particle or an approximation
(regularization) to a delta function. Rather, a blob can be thought
of as a diffuse sphere that has some physical extent and interacts
with the continuum in its interior.

The minimally-resolved blob approach can be used in a variety of contexts
where boundary conditions are specified on the surface of particles
embedded in a medium. These boundary conditions can be of the Dirichlet,
Neumann, or mixed (Robin) type. Boundary conditions involving normal
or tangential derivatives could be specified on the exterior boundary
of the particles only (e.g., free-slip for flow problems), or as a
jump condition between the interior and exterior (e.g., thermal or
electrical conductivity in a dispersion of conducting particles).
The specifics of how to express the surface boundary conditions to
a volumetric blob condition are very problem specific and need to
be carefully constructed on a case by case basis.

\subsection{Suspensions of reactive rigid spheres}

In this work we develop a blob representation for dispersions of reactive
spherical particles in either two or three dimensions. The particles
are taken as (possibly overlapping) rigid spheres (disks in two dimensions)
occupying a region $\mathcal{B}\subset\Omega$ inside a domain $\Omega$,
and we consider the general case of a reaction that is not necessarily
diffusion-limited. In the region outside of the particles $\Omega\setminus\mathcal{B}$,
the diffusion of a chemical reactant is described with the diffusion
equation for the concentration of the species $c\left(\V r,t\right)$,
\begin{equation}
\partial_{t}c=\chi\grad^{2}c+s\left(\V r,t\right)\mbox{ in }\Omega\setminus\mathcal{B},\label{eq:c_t_continuum}
\end{equation}
where $\chi$ is the diffusion coefficient and $s$ is an external
source/sink density. Different boundary conditions may be specified
on the exterior boundary of the domain $\partial\Omega$, for example,
periodic boundary conditions in the case when $\Omega$ is a torus.
Let us focus on a single reactive particle. The boundary condition
at the surface of the particle $\partial\mathcal{B}$ is a Robin condition
relating the component of the diffusive flux along the normal vector
$\V n$ (chosen here to point into the interior of the particle) with
the concentration at the surface,
\begin{equation}
\chi\left(\V n\cdot\grad c\right)=k\, c\mbox{ on }\partial\mathcal{B},\label{eq:surface_BC}
\end{equation}
where $k$ is the \emph{surface} reaction rate \cite{ReactionDiffusion_Cubic,ReactionDiffusion_Lu,ReactionDiffusion_NDL}.
The relative speed of reaction and diffusion is measured by the dimensionless
Damköhler number $\text{Da}=\left(ka\right)/\chi$ , where $a$ is
the radius of the spherical particle. Here we will follow the notation
of Lu \cite{ReactionDiffusion_Cubic}, and use the dimensionless constant
$P=\chi/(ka)=\text{Da}^{-1}$. The reaction is diffusion-limited if
$P\ll1$, and reaction-limited if $P\gg1$. In the diffusion-limited
regime $P\rightarrow0$, which is of most practical importance, the
boundary condition for (\ref{eq:c_t_continuum}) becomes $c=0$ on
$\partial\mathcal{B}$. 

In many practical situations one is interested in the steady-state
solutions, in which case one can set $\partial_{t}c=0$ in (\ref{eq:c_t_continuum}).
Quantities such as the average or total reaction rate can be obtained
by solving the steady-state reaction-diffusion problem analytically,
or numerically using first-passage Monte Carlo techniques \cite{ReactionDiffusion_Cubic,ReactionDiffusion_Lu,ReactionDiffusion_NDL,FPKMC2}
or multipole expansion techniques \cite{ReactionDiffusion_Brady,BrownianDynamics_OrderNlogN}.
At steady-state, the total reactive flux at the surface of a given
particle, $\lambda=\chi\int_{\partial\mathcal{B}}\left(\V n\cdot\grad c\right)dS$,
is finite and gives the effective source/sink strength of the particle.
For an isolated sphere of radius $a$ immersed in a reservoir of the
reactant with concentration $c_{\infty}$ in the absence of sources
($s=0$), the analytical solution (see, for example, the Appendix
of Ref. \cite{ReactionDiffusion_Cubic}) has a modified Smoluchowski
form,
\begin{equation}
c\left(r\right)=c_{\infty}\left(1-\frac{1}{1+P}\frac{a}{r}\right).\label{eq:c_r_sphere}
\end{equation}
The slow $r^{-1}$ decay of the concentration is analogous to the
slow decay of hydrodynamic interactions in suspensions of spheres,
and leads to nontrivial interactions between the particles in the
dispersion. The total rate of consumption of the reactant is 
\begin{equation}
\lambda=\frac{4\pi a\chi c_{\infty}}{\left(1+P\right)}.\label{eq:lambda_sphere}
\end{equation}

Let us consider a dispersion of $N$ spheres in a domain of volume
$V_{\Omega}$ in the case when a uniform source density of reactant
is supplied, $s\left(\V r,t\right)=\text{const}$. The steady-state
average concentration $\bar{c}=V_{\Omega}^{-1}\int_{\Omega\setminus\mathcal{B}}c\, d\V r$
will be given by the balance condition between the total consumption
rate at the surfaces of the particles and the supply in the bulk.
We can define an overall reaction rate constant $K$ via the relation
\cite{ReactionDiffusion_Cubic}
\[
K\bar{c}=N^{-1}\left[\chi\int_{\partial\mathcal{B}}\left(\V n\cdot\grad c\right)dS\right]=N^{-1}s\left(V_{\Omega}-V_{\mathcal{B}}\right),
\]
where $V_{\mathcal{B}}=N\left(4\pi a^{3}/3\right)$ is the volume
occupied by the particles. Denoting the volume or packing fraction
$\varphi=V_{\mathcal{B}}/V_{\Omega}$, we can express $K$ as
\[
K=\left(4\pi a^{3}/3\right)\frac{\left(1-\varphi\right)s}{\varphi\bar{c}}.
\]
In a very dilute dispersion of well-separated spheres, $\varphi\ll1$,
we can ignore particle interactions and use (\ref{eq:c_r_sphere})
to obtain that
\[
K_{0}=\lim_{\varphi\rightarrow0}K=\frac{4\pi a\chi}{1+P}.
\]
A dimensionless (normalized) reaction rate constant that measures
the importance of multiparticle effects can be defined as \cite{ReactionDiffusion_Cubic}
\begin{equation}
\beta_{P}=\frac{K}{K_{0}}=\frac{sa^{2}}{\chi}\cdot\frac{\left(1-\varphi\right)}{3\varphi}\cdot\frac{\left(1+P\right)}{\bar{c}},\label{eq:beta_P_cont}
\end{equation}
and measured by calculating the average concentration $\bar{c}$ at
steady state. The dimensionless quantity $\beta_{P}$ is independent
of the particular physical units and only depends on the dimensionless
number $P$ and the microstructure of the system, i.e., on $\varphi$
and the arrangement of the spheres. The value of this quantity has
been obtained using first-passage techniques for a variety of regular
(lattice) arrangements of spheres and for random dispersions \cite{ReactionDiffusion_Cubic,ReactionDiffusion_Lu,ReactionDiffusion_NDL,ReactionDiffusion_Brady,ReactionDiffusion_Torquato}.

As a model problem on which we will validate and calibrate the blob
model, let us consider the steady state solution for a finite collection
of $N$ reactive spheres of radius $a$ in an infinite three-dimensional
reservoir with concentration $c_{\infty}$ far away from the particles.
If one approximates each reactive sphere $i$ with just a monopole
term of strength $\lambda_{i}=\chi\int_{\partial\mathcal{B}_{i}}\left(\V n\cdot\grad c\right)dS$,
where $\partial\mathcal{B}_{i}$ is the surface of the sphere, then
a standard multipole expansion gives the system of equations (see
Section II.B in Ref. \cite{ReactionDiffusion_Lu})
\[
\chi c_{\infty}=\frac{\lambda_{j}}{4\pi a}+\frac{1}{4\pi}\sum_{i\neq j=1}^{N}\frac{\lambda_{i}}{r_{ij}},\quad j=1,\dots,N,
\]
where $r_{ij}=\norm{\V q_{i}-\V q_{j}}$ is the distance between particles
$i$ and $j$. In matrix form we can write this as
\begin{equation}
\chi c_{\infty}\,\M 1=\M{\mathcal{M}}\V{\lambda},\label{eq:monopole_eq}
\end{equation}
where $\M 1$ is a vector of ones and the matrix $\M{\mathcal{M}}$
has entries $\mathcal{M}_{ii}=\left(4\pi a\right)^{-1}$ and $\mathcal{M}_{ij}=\left(4\pi r_{ij}\right)^{-1}$
for $i\neq j$. The slow $r^{-1}$ decay is a characteristic signature
of the inverse Poisson operator and couples the particles through
their far field contributions. Note that in two dimensions the Green's
function for the Poisson equation behaves logarithmically and the
boundary conditions affect the result strongly regardless of the system
size. In fact, in two dimensions one cannot prescribe the concentration
at infinity (this is related to the Stokes paradox for flow around
a cylinder).

\section{Reactive Blob Model}

In the beginning, we focus on the continuum formulation of the continuum-particle
coupling. However, it is important to point out that most of the notation
and conclusions can directly be adopted in the discrete formulation
by simply replacing spatial integrals with sums over grid points.
We will return to the spatially-discrete formulation in Section \ref{sec:Discretization}.

\subsection{The Blob Model of a Particle}

Let us consider a single spherical reactive particle of radius $a$
at position $\V q(t)$ in a domain with Cartesian coordinate $\V r$.
The shape of the particle and its effective interaction with the surrounding
medium is captured through a smooth kernel function $\delta_{a}\left(\V r\right)$
that integrates to unity and whose support is localized in a region
of size $a$. In the continuum setting considered here, one may choose
any one-dimensional ``bell-shaped'' curve $\delta_{a}\left(r\right)$
with half-width of order $a$, and define a spherically-symmetric
$\delta_{a}\left(\V r\right)=\delta_{a}\left(r\right)$; alternatively,
in $d$ dimensions one may define a tensor-product
\begin{equation}
\delta_{a}\left(\V r\right)=\prod_{\alpha=1}^{d}\delta_{a}\left(r_{\alpha}\right).\label{eq:delta_tensor}
\end{equation}
In immersed-boundary methods \cite{IBM_PeskinReview}, the kernel
function $\delta_{a}$ is considered to be an approximation of the
Dirac delta function of purely numerical origin and has the tensor-product
form (\ref{eq:delta_tensor}). Here we choose the shape of the function
based on numerical considerations, but relate its shape to the physical
properties of the particle, namely, its effective reactive radius,
surface and volume. This is very similar to the approach taken for
particulate flows in the Force Coupling Method \cite{ForceCoupling_Monopole}.

The reactant concentration field $c(\V r,t)$ is extended over the
whole domain, \emph{including} the particle interior. Interaction
between the field and particle is mediated by the kernel function
through two crucial local operations. The\emph{ local averaging }operator\textbf{
$\M J(\V q)$} averages the concentration inside the particle to estimate
the local concentration
\[
c_{\V q}\left(t\right)=\int\delta_{a}\left(\V q-\V r\right)c\left(\V r,t\right)\, d\V r\equiv\left[\M J(\V q)\right]c.
\]
The reverse of local averaging is accomplished using the \emph{local
spreading }linear operator $\M S(\V q)$ which takes a point source/sink
strength $\lambda\left(t\right)$ at the location of the particle
and spreads it over the extent of the kernel function to return a
smooth source density field,
\[
\lambda_{\V q}\left(\V r,t\right)=\lambda\left(t\right)\,\delta_{a}\left(\V q-\V r\right)\equiv\left[\M S(\V q)\right]\lambda.
\]
Note that the local spreading operator $\M S$ has dimensions of inverse
volume. For notational simplicity we will slightly abuse notation
and assume that the local spreading and interpolation operators can
be applied to a scalar, a vector, or a tensor field, with the interpretation
that the same local averaging or spreading operation is applied to
each component independently. This sort of block-diagonal form of
the spreading and interpolation operators is not strictly required
for the mathematical formulation, but applies to the specific forms
of the operators we use in practice \cite{IBM_PeskinReview}.

The physical volume of the particle $\D V$ is related to the shape
and width of the kernel function via $\M J\M S=\D V^{-1}\,\M I$,
that is,
\begin{equation}
\D V=\left(\M J\M S\,1\right)^{-1}=\left[\int\delta_{a}^{2}\left(\V r\right)d\V r\right]^{-1}.\label{eq:dV_JS}
\end{equation}
Therefore, even though the particle is represented only by the position
of its centroid, it is not appropriate to consider it a ``point''
particle. Rather, it can be thought of as a diffuse sphere that has
some physical extent. For lack of better terminology, we will refer
to such a diffuse particle as a ``blob''.

For multi-particle problems, the positions of $N$ blobs is described
by the configuration $\V Q=\left\{ \V q_{1},\dots,\V q_{N}\right\} $.
We can define a composite local averaging operator $\M{\mathcal{J}}\left(\V Q\right)$
that separately averages the concentration in the neighborhood of
each particle, 
\[
\left(\M{\mathcal{J}}c\right)_{i}\equiv\left[\M J(\V q_{i})\right]c,
\]
as well as a composite local spreading operator that adds the contributions
from all the sources/sinks $\V{\lambda}=\left\{ \lambda_{1},\dots,\lambda_{N}\right\} $,
\[
\M{\mathcal{S}}\V{\lambda}=\sum_{i=1}^{N}\left[\M S(\V q_{i})\right]\lambda_{i}.
\]
It is important to note that the averaging and spreading operators
are \emph{adjoint}, $\M{\mathcal{S}}=\M{\mathcal{J}}^{\star}$. The
corresponding spatially-discrete operators (matrices) will be scaled
transposes of each other. Also note that if the kernels are compactly-supported
and the supports of the kernels of the different particles are not
overlapping, $\M{\mathcal{J}}\M{\mathcal{S}}=\D V^{-1}\,\M I_{N}$
is a multiple of the identity operator.

\subsection{\label{sub:ReactiveBlobs}Reaction-Diffusion Equation with Blobs}

Let us consider a single reactive blob. Instead of coupling the particle
and the medium via surface integrals, as in (\ref{eq:c_t_continuum},\ref{eq:surface_BC}),
we model the reaction-diffusion system by
\[
\partial_{t}c=\chi\grad^{2}c-\sum_{i}\kappa_{i}\left[\int\delta_{a}\left(\V q_{i}-\V r\right)c\left(\V r,t\right)\, d\V r\right]\delta_{a}\left(\V q_{i}-\V r\right)+s\left(\V r,t\right),
\]
where $\kappa_{i}$ is an overall reaction rate near blob $i$. Using
the local interpolation and spreading operators the reaction-diffusion
equation can be succinctly written as
\begin{equation}
\partial_{t}c=\chi\grad^{2}c-\left(\sum_{i=1}^{N}\kappa_{i}\M S_{i}\M J_{i}\right)c+s=\chi\grad^{2}c-\left(\M{\mathcal{S}}\M{\kappa}\M{\mathcal{J}}\right)c+s,\label{eq:c_t_blob}
\end{equation}
where $\M{\kappa}$ is a diagonal matrix of reaction rates. Note that
when particles consume reactant, $\kappa>0$, the operator $-\M{\mathcal{S}}\M{\kappa}\M{\mathcal{J}}$
is negative-semidefinite just like the Laplacian operator, and therefore
the dynamics (\ref{eq:c_t_blob}) is strictly dissipative.

The surface reaction rate $k$ in (\ref{eq:c_t_continuum},\ref{eq:surface_BC})
has physical units of $m/s$. The blob reaction rate $\kappa=kA\sim k\D V/a\sim ka^{2}$
has units $m^{3}/s$, where $A$ is the surface area of the blob.
In order to relate the blob model to the perfect sphere model, we
take
\begin{equation}
\kappa=kA=4\pi ka^{2},\label{eq:kappa_k}
\end{equation}
where $a$ is an effective reactive radius of the blob, related to
the width of the kernel function $\delta_{a}$. We will justify the
identification (\ref{eq:kappa_k}) in Section \ref{sub:FiniteReactionRate}.
Recall that the inverse Damköhler dimensionless number 
\begin{equation}
P=\frac{\chi}{ka}=4\pi\frac{\chi a}{\kappa}\label{eq:P_def}
\end{equation}
determines whether the process is diffusion-limited ($P\rightarrow0$)
or reaction-limited ($P\rightarrow\infty$).

In the diffusion-limited case, the total rate of consumption of the
reactant by particle $i$, $\lambda_{i}=\kappa_{i}\left(\M J_{i}c\right)$,
remains constant even though $\kappa_{i}\rightarrow\infty$. This
means that $\M J_{i}c\rightarrow0$, which corresponds to the boundary
condition $c=0$ on $\partial\mathcal{B}$ in the case of a fully-resolved
reactive sphere. In the limit $P\rightarrow0$ the reaction-diffusion
equation (\ref{eq:c_t_blob}) approaches the saddle-point problem
\begin{align}
\partial_{t}c & =\chi\grad^{2}c-\M{\mathcal{S}}\V{\lambda}+s,\nonumber \\
\mbox{s.t. }\M{\mathcal{J}}c & =0,\label{eq:saddle_point_c}
\end{align}
where the sink strengths $\V{\lambda}\leftarrow\M{\kappa}\M{\mathcal{J}}c$
are a Lagrange multiplier corresponding to the constraint. We remark
that (\ref{eq:c_t_blob}) can be viewed as a penalty method for the
constrained formulation (\ref{eq:saddle_point_c}). Even for problems
with a large but finite reaction rates, it may be preferable to use
the saddle-point formulation because (\ref{eq:c_t_blob}) becomes
increasingly ill conditioned the limit of large $\kappa$. We also
remark that effective solvers for (\ref{eq:c_t_blob}) may serve as
effective preconditioners for (\ref{eq:saddle_point_c}), although
this has not yet been attempted by us.

\section{\label{sec:Discretization}Spatio-Temporal Discretization}

In this section we describe our spatial discretization of (\ref{eq:c_t_blob})
and (\ref{eq:saddle_point_c}), and then focus on methods for efficiently
solving the resulting linear systems of equations using multigrid-based
preconditioning.

\subsection{\label{sub:SpatialDiscretization}Spatial Discretization}

Our spatial discretization of the equations (\ref{eq:c_t_blob}) is
based on standard finite-volume techniques for the Poisson equation,
combined with techniques from the immersed boundary method \cite{IBM_PeskinReview}.
In this work we focus on a uniform Cartesian grid of spacing $h$,
although it is possible to extend the approach to non-uniform adaptively-refined
grids \cite{StaggeredIBM,IBAMR}.

In the discrete setting, the various continuum operators acting on
vector fields become matrices. The concentration field $c$ is discretized
as a vector $\V c$ of cell-centered concentrations. The Laplacian
$\grad^{2}c$ is discretized as $\M L\V c$ using the standard second-order
$\left(2d+1\right)$-point discrete Laplacian $\M L$, along with
suitable modifications near physical boundaries depending on the boundary
condition imposed. Higher-order or more isotropic discrete Laplacians
can be used to lessen the discretization artifacts, but at the cost
of more complicated multigrid linear solvers and difficulties with
conservative adaptive-mesh refinement.

Application of the local averaging operator $\M J$, which is a convolution
operator in the continuum setting, becomes a discrete summation over
the grid points that are near the particle,
\[
\M J\V c\equiv\sum_{\V k\in\text{grid}}\phi_{w}\left(\V q-\V r_{k}\right)c_{k},
\]
where $\V r_{k}$ denotes the center of the cell (control volume)
with which $c_{k}$ is associated. Here $\phi_{w}\left(x\right)$
is a \emph{kernel} function with compact support of width $wh$. Note
that $\phi_{w}$ is related to the kernel function $\delta_{a}$ that
appeared in the continuum formulation, but is not in general the same
function. In particular, in the discrete setting the kernel $\phi_{w}$
needs to be carefully constructed taking into account the underlying
discrete grid. This grid is not rotationally invariant, and this complicates
the use of spherically-symmetric kernels.

\begin{figure}
\centering{}\includegraphics[width=0.45\columnwidth]{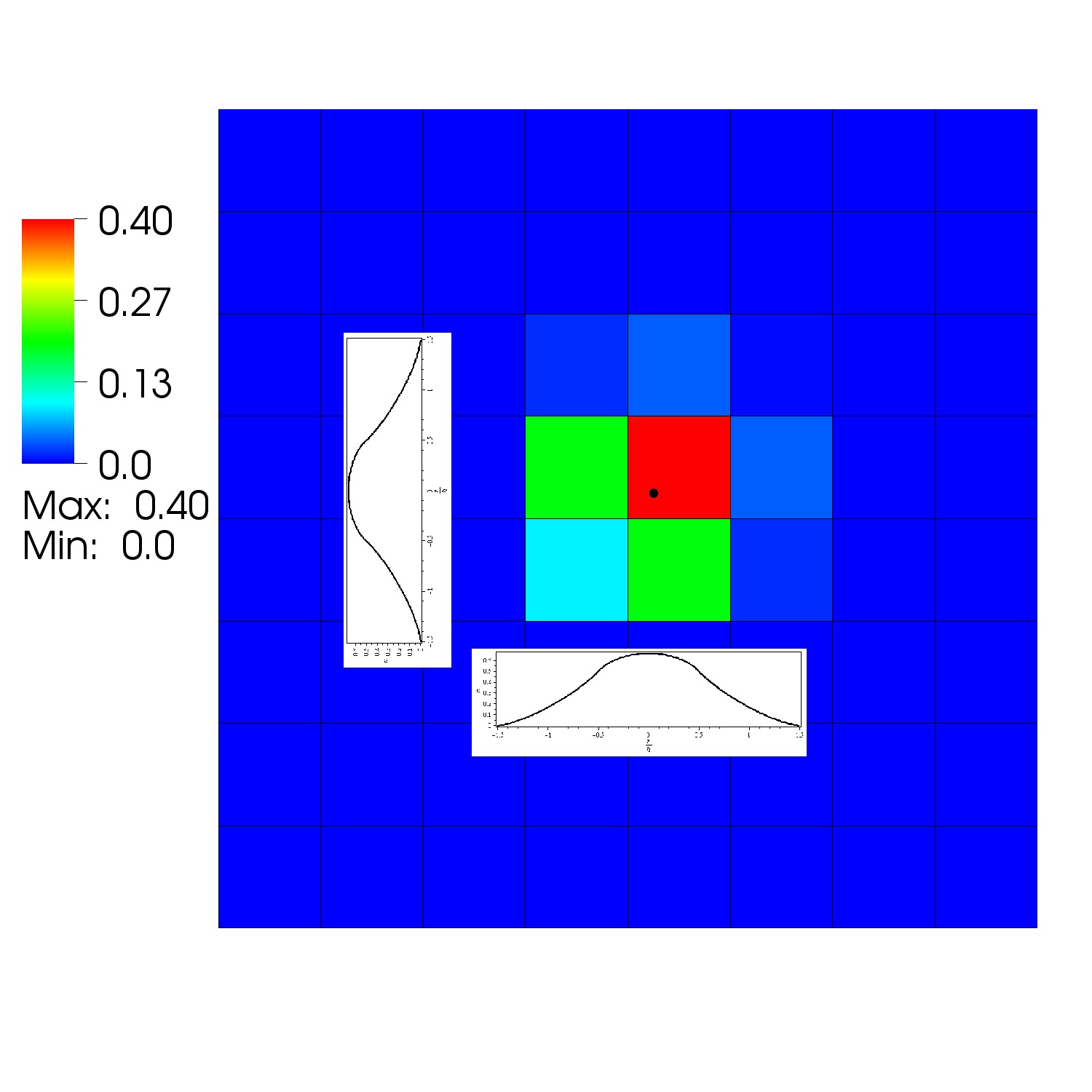}\includegraphics[width=0.45\columnwidth]{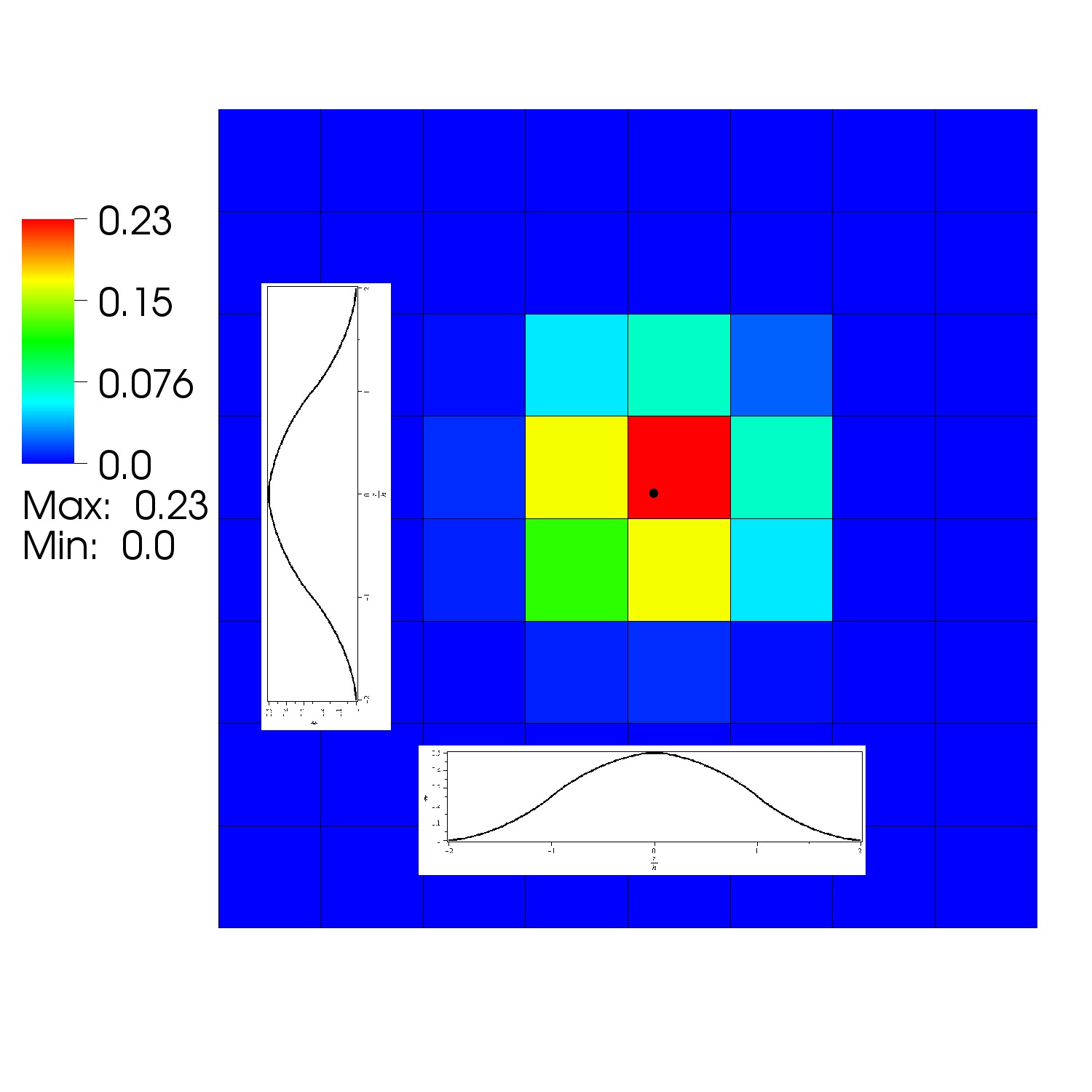}
\caption{\label{fig:illustration2D}Illustration of discrete kernel functions
used to represent the interaction between the particles and the grid
used to solve the concentration equation. The position of the blob
$\V q$ is shown with a dot. The color of each cell with center $\V r_{k}$
corresponds to the value of the kernel $\phi_{w}\left(\V q-\V r_{k}\right)$
weighting that cell. The one dimensional kernel function used to construct
the two-dimensional kernel via the tensor product (\ref{eq:tensor_product})
is shown as an inset. (\emph{Left panel}) The three-point ($w=3$)
Peskin kernel $\varphi_{3}$. (\emph{Right panel}) The four-point
($w=4$) Peskin kernel $\varphi_{4}$.}
\end{figure}

We follow the approach generally employed by immersed boundary methods
\cite{IBM_PeskinReview} and do the local averaging independently
along each of the $d$ dimensions,
\begin{equation}
\phi_{w}\left(\V q-\V r_{k}\right)=\prod_{\alpha=1}^{d}\phi_{w}\left[q_{\alpha}-\left(r_{k}\right)_{\alpha}\right],\label{eq:tensor_product}
\end{equation}
which improves the isotropy of the spatial discretization. Constructing
kernels $\phi_{w}$ that provide both translational and rotational
invariance would require not using tensor-product kernels; this difficult
task has not, to our knowledge, been accomplished yet. As a matrix,
the local spreading operator $\M S=\M J^{\star}=\left(\D V_{f}\right)^{-1}\M J^{T}$,
\[
\left(\M S\lambda\right)_{k}=\left(\D V_{f}\right)^{-1}\phi_{w}\left(\V q-\V r_{k}\right)\lambda,
\]
where $\D V_{f}=h^{d}$ is the volume of a hydrodynamic cell. Examples
of two compact-support ($w=3$ and $w=4$) discrete kernel functions
in two dimensions are illustrated in Fig. \ref{fig:illustration2D}
as color plots of the value $\phi_{w}\left(\V q-\V r_{k}\right)$
corresponding to the cell centers $\V r_{k}$ around a blob. The one
dimensional kernel function used to construct the two-dimensional
kernel via the tensor product (\ref{eq:tensor_product}) is shown
for comparison.

The kernel function $\phi_{w}$ was constructed by Peskin \cite{IBM_PeskinReview}
to yield translationally-invariant zeroth- and first-order moment
conditions, along with a quadratic condition, 
\begin{eqnarray}
\sum_{\V k\in\text{grid}}\phi_{w}\left(\V q-\V r_{k}\right) & = & 1\nonumber \\
\sum_{\V k\in\text{grid}}\left(\V q-\V r_{k}\right)\phi_{w}\left(\V q-\V r_{k}\right) & = & 0\nonumber \\
\sum_{\V k\in\text{grid}}\phi_{w}^{2}\left(\V q-\V r_{k}\right) & = & \D V^{-1}=\mbox{const.},\label{eq:JS_invariance}
\end{eqnarray}
independent of the position of the particle $\V q$ relative to the
underlying (fixed) velocity grid. Ensuring these properties requires
making the support of the kernel function an integer multiple $w$
of the grid spacing $h$, specifically, taking $\phi_{w}\left(x\right)=\tilde{\phi}_{w}\left(x/h\right)$,
where $\tilde{\phi}_{w}$ is compactly supported on the interval $\left[-w/2,\, w/2\right]$.
This means that the size and shape of the particle $a\sim h$ is determined
by the spatial discretization of the diffusion equation and the kernel
used to transfer information between the particle and the grid, and
cannot be chosen arbitrarily. We will obtain the precise relationship
between the reactive radius $a$ of a single reactive blob and the
grid spacing $h$ in Section \ref{sec:Results}.

The last condition (\ref{eq:JS_invariance}), was imposed by Peskin
\cite{IBM_PeskinReview} as a way of approximating independence under
shifts of order of the grid spacing. This property implies that the
particle volume $\D V=\left(\M J\M S\,1\right)^{-1}$ will remain
constant and independent of the position of the blob relative to the
underlying grid. The function with minimal support that satisfies
(\ref{eq:JS_invariance}) is uniquely determined \cite{IBM_PeskinReview}.
In some of our numerical experiments we employ this \emph{three-point}
discrete kernel function $\varphi_{3}$ \cite{StaggeredIBM,DirectForcing_Balboa},
illustrated in the left panel of Fig. \ref{fig:illustration2D}. This
particular choice gives $\D V=2^{d}\D V_{f}=2^{d}h^{d}$, where $d$
is the dimensionality. We find, consistent with Ref. \cite{ISIBM},
better translation when using the Peskin \emph{four-point} discrete
kernel function $\varphi_{4}$ \cite{IBM_PeskinReview}, which gives
$\D V=\left(8/3\right)^{d}h^{d}$, and is illustrated in the right
panel of Fig. \ref{fig:illustration2D}.

\subsection{Temporal Discretization}

In this section, we describe how to integrate the spatially-discretized
equations from time $n\D t$ to time $\left(n+1\right)\D t$, where
$\D t$ is the time step size. We will use a superscript to denote
the time level at which a given quantity is evaluated.

As a typical temporal integrator for (\ref{eq:c_t_blob}), let us
consider using a backward Euler step for $\V c$,
\begin{equation}
\frac{\V c^{n+1}-\V c^{n}}{\D t}=\chi\M L\V c^{n+1}-\M{\mathcal{S}}^{n}\M{\kappa}\M{\mathcal{J}}^{n}\V c^{n+1}+\V s^{n}.\label{eq:c_np1}
\end{equation}
The implicit update (\ref{eq:c_np1}) requires solving a linear system
of the form
\begin{equation}
\left(\D t^{-1}\M I-\chi\M L+\M{\mathcal{S}}\M{\kappa}\M{\mathcal{J}}\right)\V c=\M B\V c=\V g,\label{eq:unsteady_system}
\end{equation}
which we discuss in the next section. The relative importance of the
Laplacian term and the scaled identity term is measured by the diffusive
CFL number $\vartheta=\chi\D t/h^{2}$. The Euler scheme is only first-order
accurate. Second-order accuracy can be achieved by using the implicit
midpoint rule for the diffusive term, however, the advantage of using
the backward Euler method is that one can take very large time steps,
$\D t\rightarrow\infty$. Namely, in the limit $\vartheta\rightarrow\infty$
the backward Euler method approaches a steady-state solver,
\begin{equation}
\left(\M{\mathcal{S}}\M{\kappa}\M{\mathcal{J}}-\chi\M L\right)\V c=\V s.\label{eq:steady_state}
\end{equation}

\subsection{Diffusion-Limited Case}

In the limit $P\rightarrow0$, we need to solve a saddle-point problem
(\ref{eq:saddle_point_c}) for $\V c$ and $\V{\lambda}$,
\begin{align}
\M A\V c+\zeta\M{\mathcal{S}}\V{\lambda} & =\V g\nonumber \\
\xi\M{\mathcal{J}}\V c & =\V f,\label{eq:saddle_point}
\end{align}
where $\V{\lambda}\leftarrow\zeta^{-1}\M{\kappa}\M{\mathcal{J}}\V c$
is the uknown sink strength (rate of consumption of reactant) for
each of the $N$ blobs. Here $\M A=\D t^{-1}\M I-\chi\M L$, and we
have allowed for a nonzero $\V f$ since this generality is necessary
when constructing preconditioners for solving the saddle-point problem.
We have added two scaling constants $\zeta$ and $\xi=\zeta/\D{V_{f}}$
to make the system symmetric and potentially improve its numerical
scaling in cases when the different physical values have vastly different
units (a reasonable choice is $\zeta\sim\chi h$), as an alternative
to non-dimensionalization of the equations. Notice that if we impose
periodic or homogeneous Neumann boundary conditions in the steady-state
case, $\M A=-\chi\M L$, the matrix $\M A$ has a null space consisting
of constant vectors, and the range of $\M A$ is the set of vectors
of mean zero. Let us first consider the case when $\D t$ is finite
and $\M A$ is invertible.

To construct a preconditioner for the saddle-point problem, we write
down the exact solution using a Schur complement approach. From the
first equation,
\begin{equation}
\V c=\M A^{-1}\left(\V g-\zeta\M{\mathcal{S}}\V{\lambda}\right),\label{eq:c_exact}
\end{equation}
and plugging this into the constraint we get
\[
\xi\zeta\left(\M{\mathcal{J}}\M A^{-1}\M{\mathcal{S}}\right)\V{\lambda}=\xi\M{\mathcal{J}}\M A^{-1}\V g-\V f=\V h.
\]
This gives the solution for the Lagrange multipliers $\V{\lambda}$
in terms of the inverse of the Schur complement $\M M=\M{\mathcal{J}}\M A^{-1}\M{\mathcal{S}}$,
\begin{equation}
\V{\lambda}=\left(\zeta\xi\right)^{-1}(\M{\mathcal{J}}\M A^{-1}\M{\mathcal{S}})^{-1}\V h=\left(\zeta\xi\right)^{-1}\M M^{-1}\V h.\label{eq:r_schur}
\end{equation}

To understand the physical meaning of the matrix $\M M$, consider
a finite collection of $N$ reactive blobs in an infinite three-dimensional
reservoir with concentration $c_{\infty}$ far away from the particles.
Setting $\zeta=\xi=1$ and considering steady state $\M A=-\chi\M L$,
the excess concentration $\d{\V c}=\V c-c_{\infty}$ solves the saddle-point
problem 
\begin{align}
\chi\M L\left(\d{\V c}\right) & =\M{\mathcal{S}}\V{\lambda},\nonumber \\
\M{\mathcal{J}}\left(\d{\V c}\right) & =-c_{\infty},\label{eq:saddle_point_c-1}
\end{align}
the solution of which can be read from (\ref{eq:r_schur}) to be 
\[
-\left(\M{\mathcal{J}}\M L^{-1}\M{\mathcal{S}}\right)\V{\lambda}=\M M\V{\lambda}=\chi c_{\infty}.
\]
This equation can be seen as a discretization of the monopole expansion
(\ref{eq:monopole_eq}), with $\M M$ being a particular discretization
of $\M{\mathcal{M}}$. The diagonal values will be $M_{ii}\sim h^{-1}$,
to within finite-size corrections. In fact, the value of the diagonal
element of $\M M$ can be used to define an effective reactive radius
$a\sim h$ for a blob, as we will study in more detail in Section
\ref{sec:Results}. For two blobs $i$ and $j$ separated by a distance
$r_{ij}\gg h$, we have that $M_{ij}\sim r^{-1}$ just like $\mathcal{M}_{ij}$.
When the two blobs are close to each other, the $r^{-1}$ singularity
of the Green's function for the Poisson equation is regularized by
the shape of the kernel used to describe the blob. This is very similar
to how the method of regularized Stokeslets \cite{RegularizedStokeslets}
regularizes the singular Green's function for Stokes flow, and how
the Rotne-Prager tensor or similar quasi-Gaussian regularizations
are used to regularize the mobility in Brownian dynamics simulations
\cite{BrownianDynamics_OrderN}. Note that in two dimensions the decay
is logarithmic in $r$ and thus much slower than $r^{-1}$, making
the influence of the boundaries felt throughout the whole system.

\begin{figure}
\centering{}\includegraphics[width=0.45\columnwidth]{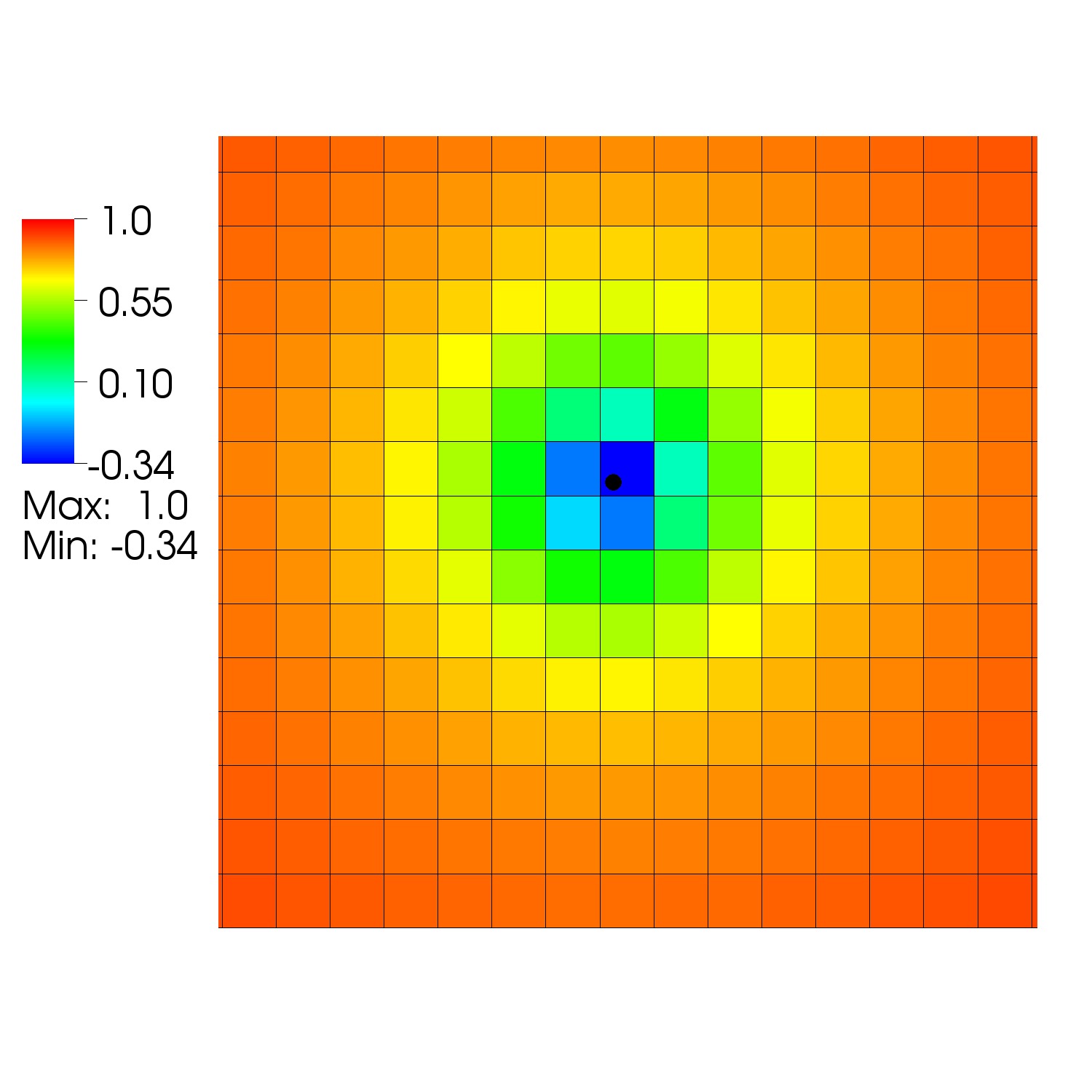}\includegraphics[width=0.45\columnwidth]{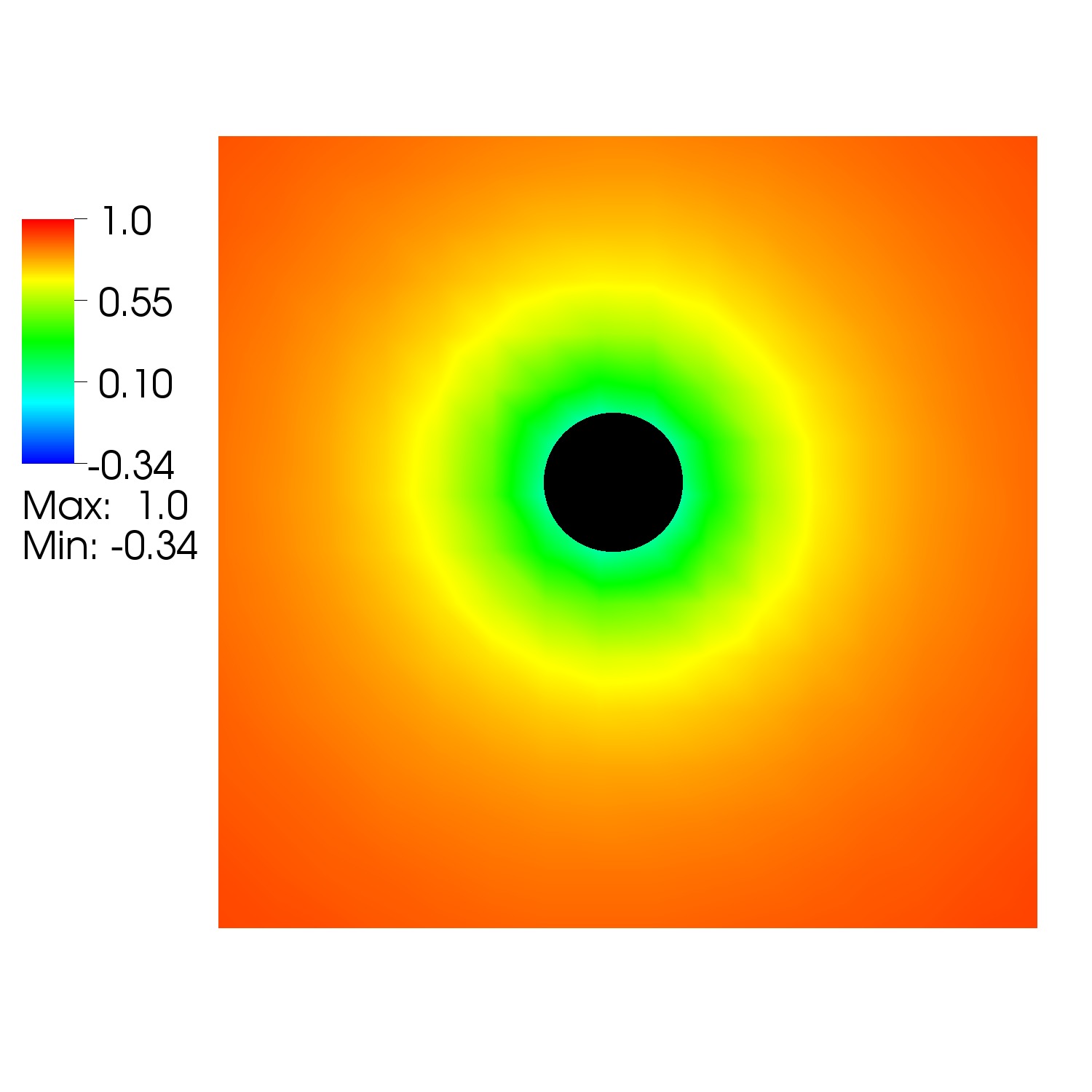}
\caption{\label{fig:Illustration3D}Ratio $c(\V r)/c_{\infty}$ for a single
blob placed at the center a large cubic box with the concentration
at the boundaries fixed to $c_{\infty}$. A slice through the three-dimensional
domain is shown. The four-point kernel $\phi_{4}$ is used as the
kernel. (\emph{Left panel}) The discrete concentration field $\V c$
(color plot) and the center of the blob (black dot). (\emph{Right
panel}) An interpolated concentration $c(\V r)$ and the equivalent
rigid reactive sphere (black disk) showing the effective reactive
radius of the blob.}
\end{figure}

As an illustration, let us consider the problem of a single reactive
sphere in an infinite three-dimensional reservoir of reactant keeping
the concentration at infinity fixed at $c_{\infty}$. The steady state
solution for the concentration is given in (\ref{eq:c_r_sphere})
with $P=0$ and has a characteristic $r^{-1}$ decay away from the
surface of the sphere. We mimic this situation by placing a single
blob near the center a large cubic box with the concentration at the
boundaries fixed to $c_{\infty}$, and set $\V f=0$. The solution
of the steady-state ($\M A=-\chi\M L$) diffusion-limited problem
(\ref{eq:saddle_point}) is illustrated in Fig. \ref{fig:Illustration3D}.
A two-dimensional slice containing the center of the blob is shown.
The panel on the left shows the discrete concentration field $\V c$,
while the the panel on the right shows an interpolated concentration
$c(\V r)$ corresponding to the discrete $\V c$. It can be seen in
the left panel that some of the concentrations near the center of
the blob are negative, which is clearly unphysical. That this has
to be so is easily seen from the fact that $\M J\V c=0$ and all weights
in the definition of the local averaging operator $\M J$ are positive.
In the right panel we show the equivalent reactive rigid sphere as
a black disk (the procedure used to calculate its radius $a$ is detailed
in Section \ref{sub:ReactionRadius}), and it is seen that the concentration
outside of the blob kernel is positive, and in fact, we will show
in Section \ref{sub:ReactionRadius} that the numerical $c(\V r)$
quantitatively matches the field around a reactive rigid sphere. This
illustrates the fact that while the minimally-resolved blob approach
does not resolve the details in the very vicinity of the particles,
it reproduces the correct far field response.

\subsection{\label{sub:FiniteReactionRate}Finite Reaction Case}

To relate the parameters $k$ and $\kappa$, consider the case of
an isolated particle in a large periodic domain with a constant supply
of reactant $s=V_{\Omega}^{-1}$. For a very large domain $s\rightarrow0$
and the steady state solution will be the same as for the case of
an isolated particle in an infinite reservoir with concentration at
infinity fixed at $c_{\infty}$, with the value of $c_{\infty}$ determined
from the balance of the consumption of reactant by the blob and the
total inflow of reactant $sV_{\Omega}=\lambda=1$. The condition $\lambda=1$
together with (\ref{eq:lambda_sphere}) gives
\[
c_{\infty}(P)=\frac{1}{4\pi a\chi}+\frac{P}{4\pi a\chi}=c_{\infty}\left(P=0\right)+\frac{P}{4\pi a\chi}.
\]
For a large system the average value of concentration in the domain
is the same as the far-field value $c_{\infty}$, 
\begin{equation}
\bar{c}_{P}=\bar{c}_{0}+\frac{P}{4\pi\chi a},\label{eq:c_bar_P}
\end{equation}
where $\bar{c}_{P}$ denotes the average steady-state concentration
at finite $P$, and $\bar{c}_{0}$ is the corresponding value for
$P=0$. It is this relation that can be taken as the definition of
$P$ for a blob, and (\ref{eq:kappa_k}) can be derived from this
definition, as we show next.

For a blob in a periodic domain with a constant source, $c_{P}$ is
the solution to
\[
\chi\M Lc_{P}-\kappa\M S\M Jc_{P}+s=0,
\]
while $c_{0}$ is the solution to the saddle-point problem
\[
\chi\M Lc_{0}-\M S\lambda+s=0,\quad\M Jc_{0}=0,
\]
with solution $\lambda=sV_{\Omega}=1$ determined from the overall
balance of reactant inflow and outflow. The difference $\tilde{c}=c_{P}-c_{0}$
is the solution to
\[
\chi\M L\tilde{c}=\M S\left(\kappa\M J\tilde{c}-\lambda\right).
\]
For this equation to be solvable, the right hand side must have mean
zero, which in this case of a single particle implies 
\[
\kappa\M J\tilde{c}=\lambda,\text{ and }\M L\tilde{c}=0.
\]
This means that $\tilde{c}$ is a constant, $\tilde{c}=\lambda/\kappa=\kappa^{-1}$,
which is consistent with (\ref{eq:c_r_sphere}), which gives $\tilde{c}(r)=c(r;P)-c(r;P=0)=P/\left(4\pi\chi a\right)$,
if we identify $\kappa^{-1}=P/\left(4\pi\chi a\right)$. From $\bar{c}_{P}=\bar{c}_{0}+\kappa^{-1}$
and (\ref{eq:c_bar_P}) we also get (\ref{eq:P_def}), confirming
that (\ref{eq:kappa_k}) is the correct way to relate $k$ for a rigid
sphere and $\kappa$ for a blob.

\subsection{\label{sub:Iterative-Linear-Solvers}Iterative Linear Solvers}

The most challenging aspect of the numerical algorithm is the solution
of linear systems such as (\ref{eq:unsteady_system}), (\ref{eq:steady_state}),
and (\ref{eq:saddle_point}). We solve these linear systems using
an iterative (matrix-free) Krylov solver. Because of the inherent
ill-conditioning in the equations, efficient solution by iterative
methods requires constructing good preconditioners (approximate solvers)
for the equations. In this section we discuss and test some ideas
for constructing effective preconditioners, especially in the most
difficult case of diffusion-limited steady state problems. Although
this section is somewhat technical, the construction of effective
solvers is necessary to treat multi-particle systems.

The solution of a Helmholtz (Poisson in the steady-state case) equation
\begin{equation}
\left(\D t^{-1}\M I-\chi\M L\right)\V c=\M A\V c=\V g\label{eq:Helmholtz}
\end{equation}
can be accomplished iteratively very efficiently by using standard
geometric multigrid methods. We rely on the IBAMR library \cite{IBAMR}
to provide these multigrid solvers. The multigrid method consists
of repeated application of multigrid V cycles, and typically each
V cycle reduces the norm of the residual by a factor of about $5-10$,
depending on the boundary conditions and details of the multigrid
method. An alternative method to solve the Helmholtz equation is to
use an Krylov method such as preconditioned conjugate gradients (PCG).
One or several cycles of multigrid can be used as a preconditioner.
For the simple constant-coefficient Helmholtz or Poisson problem the
overall performance of the Krylov method is similar to that of pure
multigrid, however, in many cases the Krylov iteration makes the linear
solver more robust. Since the computational cost is approximately
proportional to the total number of multigrid cycles, we use this
as a proxy for the CPU effort needed to obtain a given reduction of
the residual. In all convergence plots reported here the $y$ axis
is the relative residual (on a logarithmic scale) and the $x$ axis
is the total number of Helmholtz or Poisson multigrid cycles.

\subsubsection{Finite Reaction Rate}

To solve (\ref{eq:unsteady_system}) or (\ref{eq:steady_state}),
which are of the form $\M B\V c=\V g$, we use a Krylov method preconditioned
with $n$ cycles of multigrid for the Helmholtz problem (\ref{eq:Helmholtz}),
\begin{equation}
\M P_{n}^{-1}\approx\left(\D t^{-1}\M I-\chi\M L\right)^{-1}.\label{eq:Pinv_Poisson}
\end{equation}
Note that the preconditioned linear operator $\M P_{n}^{-1}\M B$
is not symmetric. Here we employ the left-preconditioned FGMRES method
as a robust Krylov solver. The FGMRES solver uses the modified Gram-Schmidt
orthogonalization to avoid the potential occurrence of a degenerate
Krylov basis. Note that we have tried other preconditioners, such
as a multiplicative preconditioner that separately solves the pure
diffusive and pure reactive sub-problems in sequence, but have found
no substantial improvement over the simple preconditioner (\ref{eq:Pinv_Poisson}).

\begin{figure}
\centering{}\includegraphics[width=0.5\paperwidth]{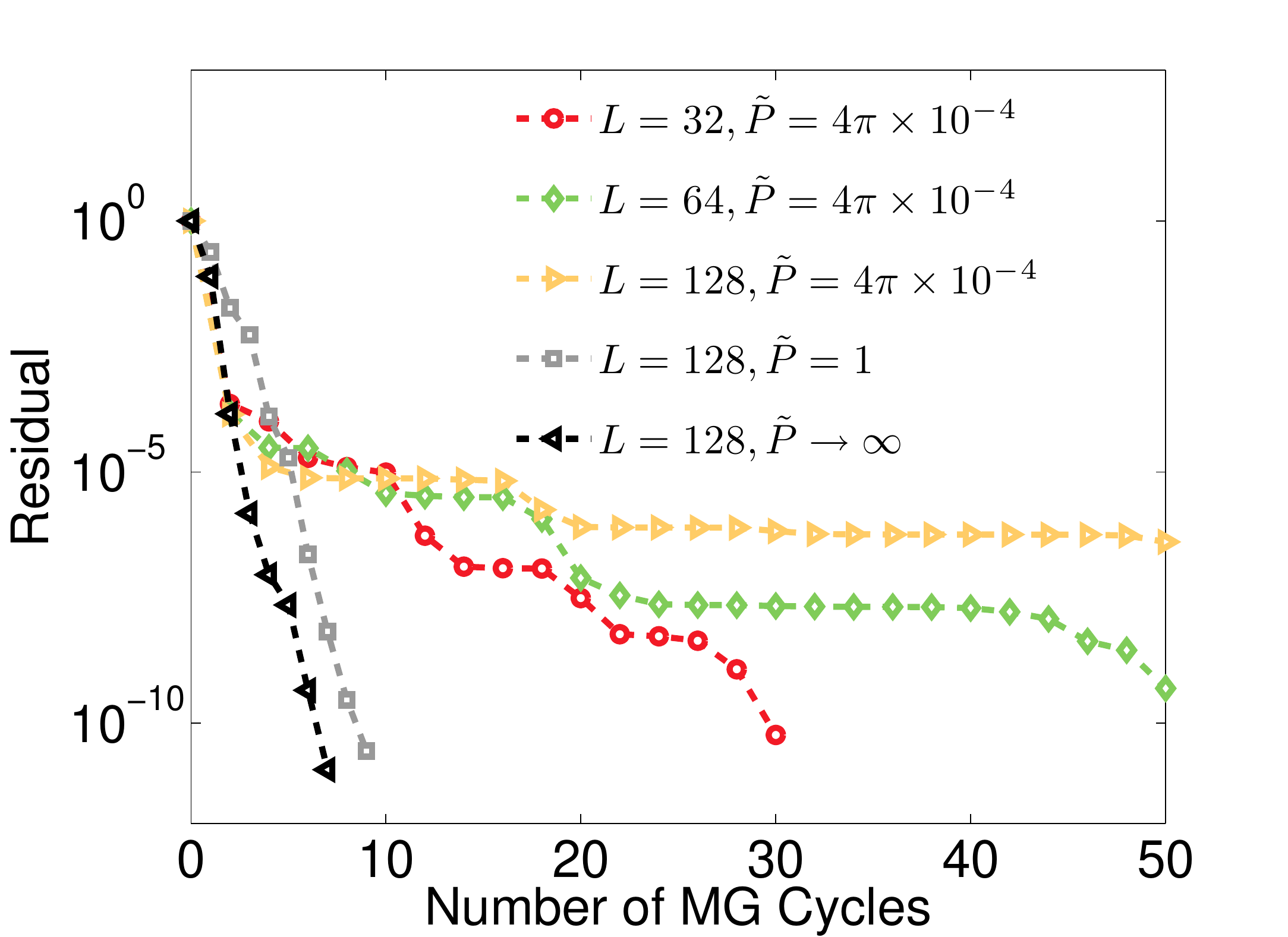}
\caption{\label{fig:NonSaddleGMRES}Convergence history of the iterative linear
solver for (\ref{eq:steady_state}) for a cubic array of $\left(L/4\right)^{3}$
blobs in a periodic domain of $L^{3}$ cells. A single V cycle of
multigrid for the Poisson equation is used as a preconditioner, and
a restart frequency of 30 is used in the FGMRES algorithm. The relative
norm of the residual versus the total number of multigrid cycles (a
proxy for the computing effort) is shown.}
\end{figure}

We study the performance of this preconditioner on the steady state
equation (\ref{eq:steady_state}) for a periodic cubic array of reactive
blobs. The periodic cubic domain is discretized into $L^{3}$ cells,
and blobs are placed 4 cells apart along each of the coordinate axes,
making the total number of blobs $\left(L/4\right)^{3}$. For these
tests we employ the 4-point kernel, so that the kernels of the blobs
are just touching but not overlapping. Let us define the dimensionless
number $\widetilde{P}=\left(4\pi\right)\chi h/\kappa\sim P$ assuming
that $a\approx h$ in (\ref{eq:P_def}). In Fig. \ref{fig:NonSaddleGMRES}
we report the performance of the multigrid preconditioner (\ref{eq:Pinv_Poisson})
for $\widetilde{P}=\left(4\pi\right)\cdot10^{-4}$, which is in the
diffusion-limited regime, for several system sizes $L$ ranging from
32 to 128 grid cells. We show the relative residual of the linear
solver as a function of the number of multigrid V cycles, which is
equivalent to the number of FGMRES iterations since we use a single
V cycle as a preconditioner, $n=1$. For comparison, in Fig. \ref{fig:NonSaddleGMRES}
we show the performance of the iterative solver for the case of $L=128$
and a reaction-limited regime $\kappa=0$, $\tilde{P}\rightarrow\infty$.
In this case the linear system (\ref{eq:steady_state}) is a simple
Poisson problem and multigrid is known to be an excellent iterative
solver. We find similarly good performance for $\tilde{P}=1$, as
also shown in the figure.

Figure \ref{fig:NonSaddleGMRES} shows that the multigrid preconditioner
is only effective for smaller system sizes or for non-diffusion-limited
systems, $\tilde{P}\gtrsim1$. For $\tilde{P}\ll1$ the solver shows
slow convergence for $L=128$, and we have not yet found any effective
simple preconditioner. Recently, specialized multigrid solvers have
been proposed for solving (\ref{eq:steady_state}) \cite{IBMultigrid_Guy,GeometricMultigrid_Boyce}.
In these methods, a matrix representation of the reaction operator
$\M{\mathcal{S}}\M{\kappa}\M{\mathcal{J}}$ is formed, and coarsened
versions of this operator are constructed using standard geometric
multigrid techniques. These methods show excellent promise and in
future work we will explore their application to reaction-diffusion
problems. Because these solution techniques \cite{IBMultigrid_Guy,GeometricMultigrid_Boyce}
consider the case of finite $\kappa$, it is not clear whether they
will aid in the solution of the diffusion-limited case $\kappa\rightarrow\infty$,
in which case it is appropriate to use a saddle-point formulation.

\subsubsection{Diffusion-Limited Case}

In the diffusion-limited case $P=0$, we need to solve the saddle-point
system (\ref{eq:saddle_point}). Formally, the solution can be found
from the inverse of the Schur complement $\M M=\M{\mathcal{J}}\M A^{-1}\M{\mathcal{S}}$
via (\ref{eq:c_exact},\ref{eq:r_schur}). In the steady-state case,
$\M M$ is a discretization and regularization of the monopole matrix
$\M{\mathcal{M}}$.

In multipole expansion methods, the $N\times N$ dense matrix $\M{\mathcal{M}}$
would be constructed (though not necessarily assembled), and then
a linear system such as (\ref{eq:monopole_eq}) solved. Sophisticated
fast-multipole methods can be used to obtain the action of $\M{\mathcal{M}}$
with cost $O\left(N\log N\right)$ or even $O\left(N\right)$ \cite{BrownianDynamics_OrderNlogN,ReactionDiffusion_Brady,BrownianDynamics_OrderN}.
With different boundary conditions, the Green's function for the Poisson
problem would have to be obtained analytically or by using boundary-integral
methods \cite{FastPoisson_FFT}. In the blob approach, the Green's
function is approximated with $\M L^{-1}$ by solving a discretization
of the Poisson equation on a grid. This is simpler to implement and
allows for much larger flexibility in handling boundary conditions,
however, it requires the development of an approximation of the action
of $\M M^{-1}$ in order to construct an effective preconditioner.
In recent years fast direct solvers, related to fast-multipole methods,
have been developed to efficiently calculate the action of $\M{\mathcal{M}}^{-1}$,
however, these methods are quite complex and still an active area
of research. Once developed such methods could be used to approximate
the action of $\M M^{-1}$.

Here we take a different approach. We use a Krylov solver to solve
the block saddle-point system (\ref{eq:saddle_point}), and use an
approximation to the inverse of the Schur complement $\M M^{-1}$
to construct a preconditioner for this Krylov solver. The simplest
possible approximation to the Schur complement matrix is to use a
diagonal preconditioner based on a single-blob approximation. Denote
the \emph{scalar }value 
\[
\gamma=-\left(\M J\M L^{-1}\M S\right)\,1,
\]
which depends on the dimensionality and the system size, and, to some
extent, on the position of the blob relative to the underlying grid.
In three dimensions $\gamma\approx\left(4\pi a\right)^{-1}$ for large
systems, where $a\sim h$ is an effective reactive radius of a blob
that we calculate in Section \ref{sec:Results}. A diagonal approximation
to the Schur complement, 
\begin{equation}
\M M^{-1}\approx\chi\gamma^{-1}\M I,\label{eq:Schur_diagonal}
\end{equation}
can be used together with (\ref{eq:r_schur}) to approximate
\[
\V{\lambda}\approx\chi\left(\zeta\xi\gamma\right)^{-1}\V h.
\]
We can then obtain an approximation to $\V c$ from (\ref{eq:c_exact})
using an inexact computation of $\M A^{-1}$. Specifically, we approximate
$\M A^{-1}\approx\widetilde{\M A}_{n}^{-1}$ by using $n$ cycles
of geometric multigrid, both when computing $\V h$, and when computing
$\V c$. We will refer to the preconditioner obtained in this manner
as the \emph{diagonal preconditioner} even though the preconditioner
matrix itself is not diagonal. Note that each Krylov iteration of
the diagonal preconditioner requires $2n$ multigrid cycles. A more
detailed description of the preconditioner steps is given in Appendix
\ref{sub:AppendixPrecond}.

Another approach is to approximate $\M M^{-1}\V h$ with an inexact
Krylov solver for solving $\M M\V{\lambda}=\V h$. This \emph{inner}
Krylov solver estimates $\V{\lambda}$, which is in turn used to estimate
$\V c$ from (\ref{eq:c_exact}) using an inexact computation of $\M A^{-1}$.
These approximate values are used to construct a preconditioner for
the outer Krylov solver, which solves the saddle-point system (\ref{eq:saddle_point}).
This \emph{approximate Schur complement preconditioner} has two integer
parameters, $n$ and $m$. The number of multigrid cycles in the approximate
Poisson solves in the outer Krylov iteration is $2n$, and the number
of cycles in the inner Krylov solver is $m$, making the total number
of multigrid cycles per outer Krylov iteration $m+2n$. A more detailed
description of the preconditioner steps is given in Appendix \ref{sub:AppendixPrecond}.

We study the performance of the diagonal and approximate Schur complement
preconditioners on the saddle-point problem (\ref{eq:saddle_point})
for a periodic cubic array of reactive blobs. The periodic cubic domain
is discretized into $L^{3}$ cells, and blobs are placed 4 cells apart
along each of the coordinate axes, making the total number of blobs
$\left(L/4\right)^{3}$. For these tests we employ the 4-point kernel,
so that the kernels of the different blobs are just touching but not
overlapping. We have performed detailed investigations of how the
number of multigrid cycles employed in the approximate Poisson solves
$m$ and $n$ affect the performance of the preconditioners; here
we summarize our main observations.

\begin{figure}
\includegraphics[width=0.49\textwidth]{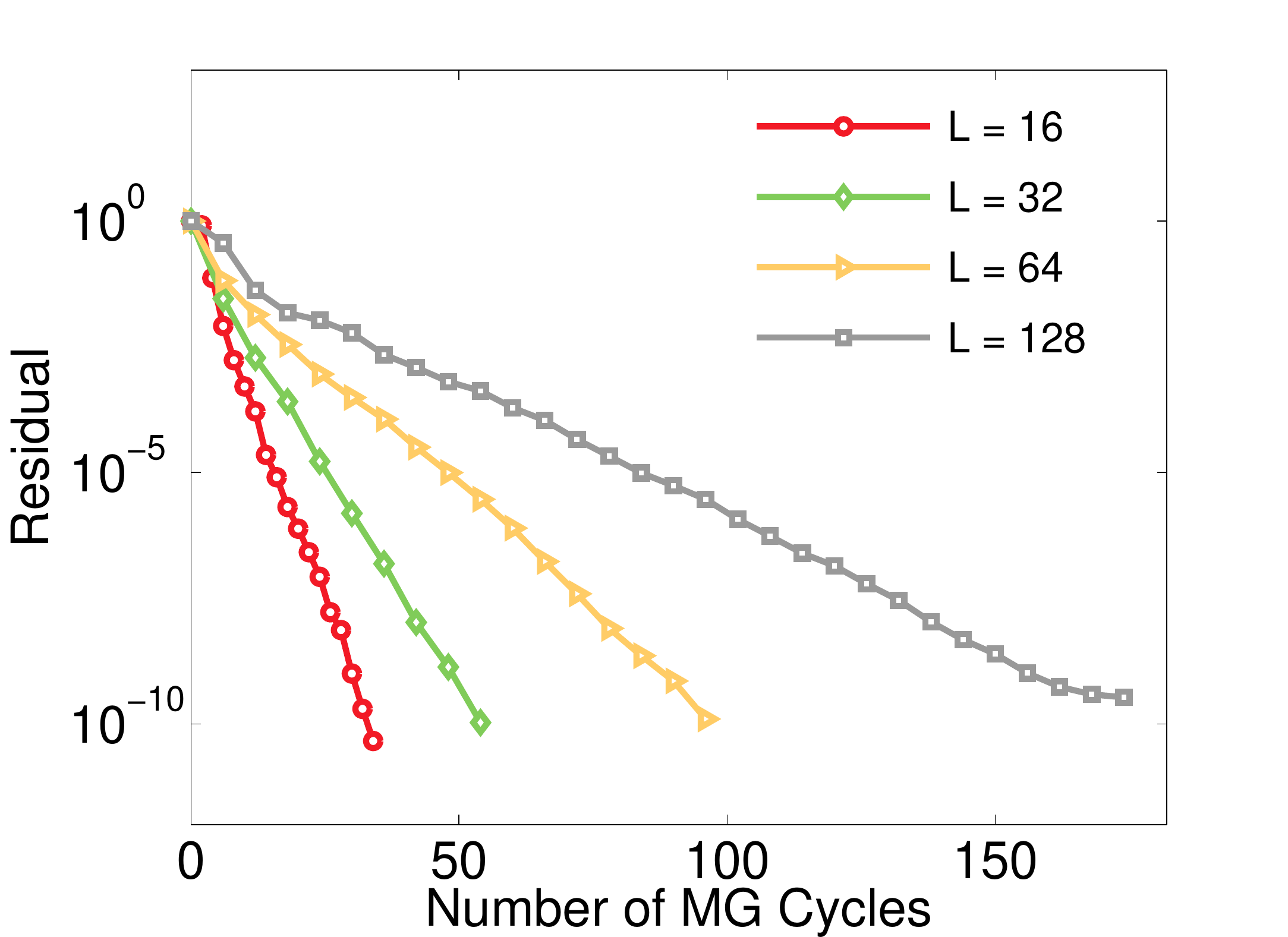}\includegraphics[width=0.49\textwidth]{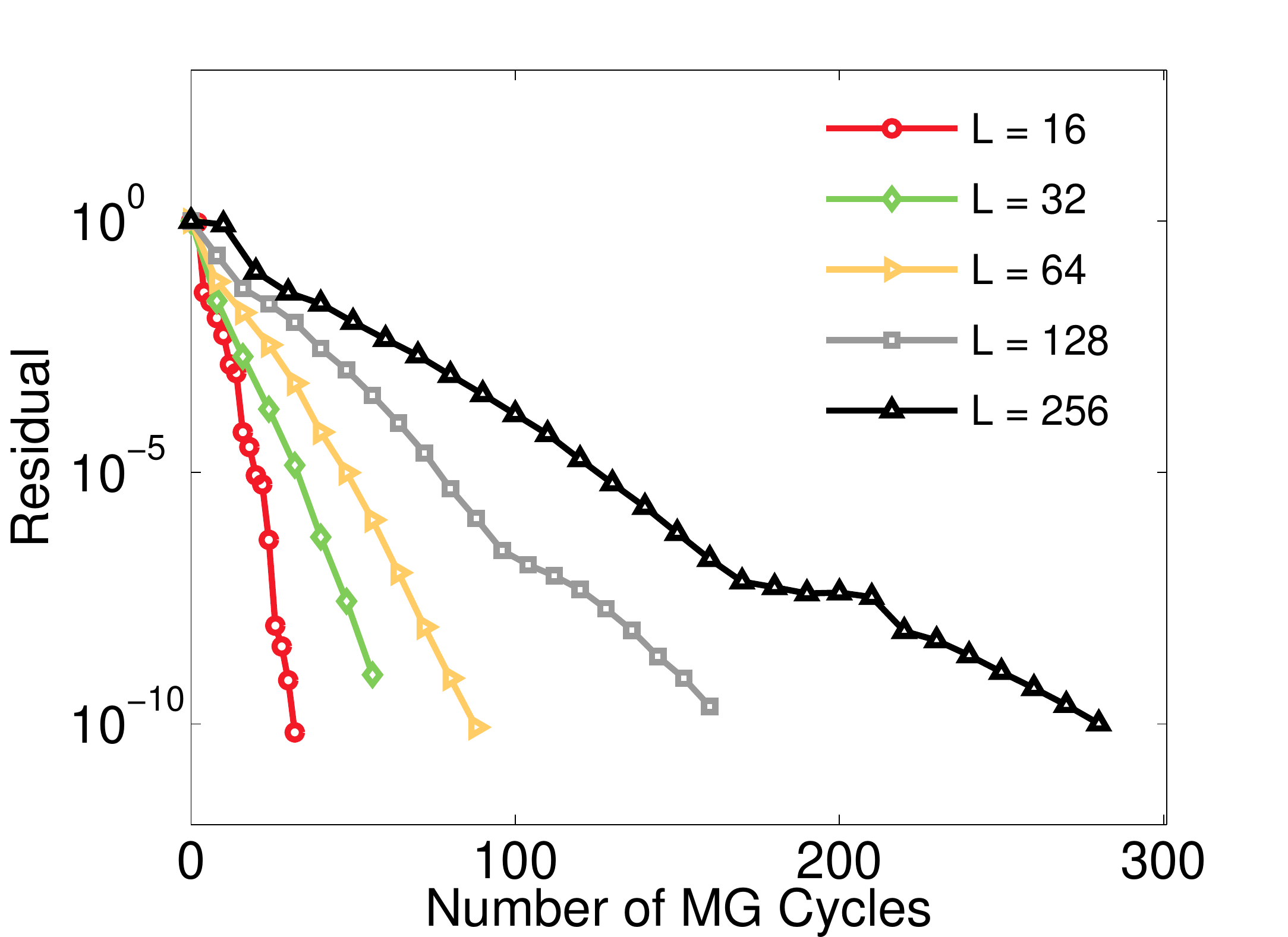}
\caption{\label{fig:DiagonalScaling}Performance of the diagonal preconditioner
with $n=1$ for the saddle-point problem (\ref{eq:saddle_point})
with $\M A=L^{-2}\M I-\M L$, as a function of the number of cells
$L$ along each dimension, in both three dimensions (left panel) and
two dimensions (right panel). The restart frequency is set to set
to 200 in the FGMRES solver.}
\end{figure}

For the first test we take $\M A=L^{-2}\M I-\M L$, which corresponds
to a very large diffusive CFL number, approaching the steady-state
problem. Note that adding a small multiple of the identity to the
Laplacian makes $\M A$ invertible even with periodic boundary conditions.
For the diagonal preconditioner we have found that a single multigrid
cycle is optimal, $n=1$ (data not shown). We have also observed that
a large restart frequency is required in the FGMRES iteration, making
the memory requirements substantial. The convergence history of the
FGMRES solver is shown in Fig. \ref{fig:DiagonalScaling} as a function
of the system size, for both two and three dimensions.

\begin{figure}
\includegraphics[width=0.49\textwidth]{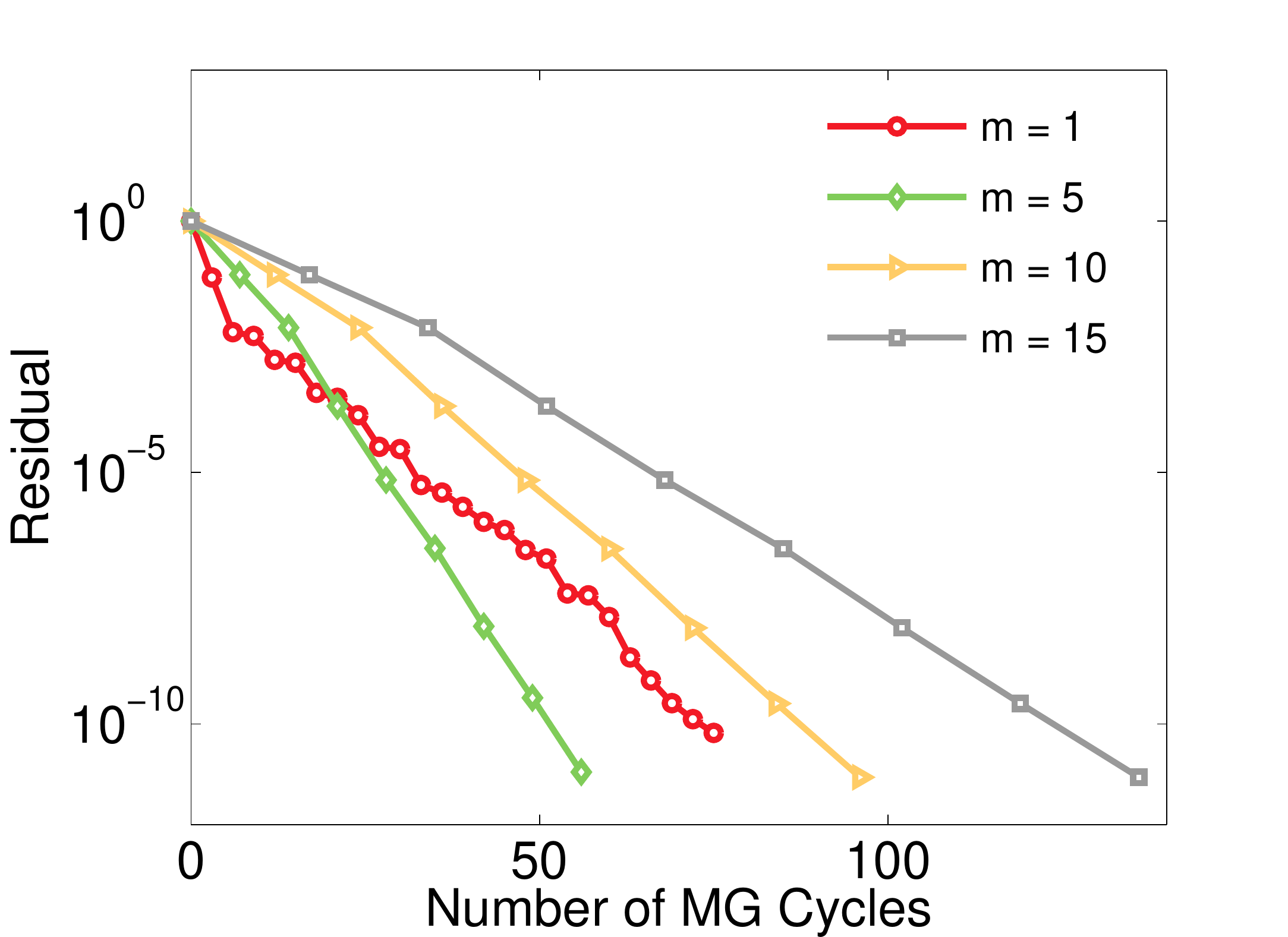}\includegraphics[width=0.49\textwidth]{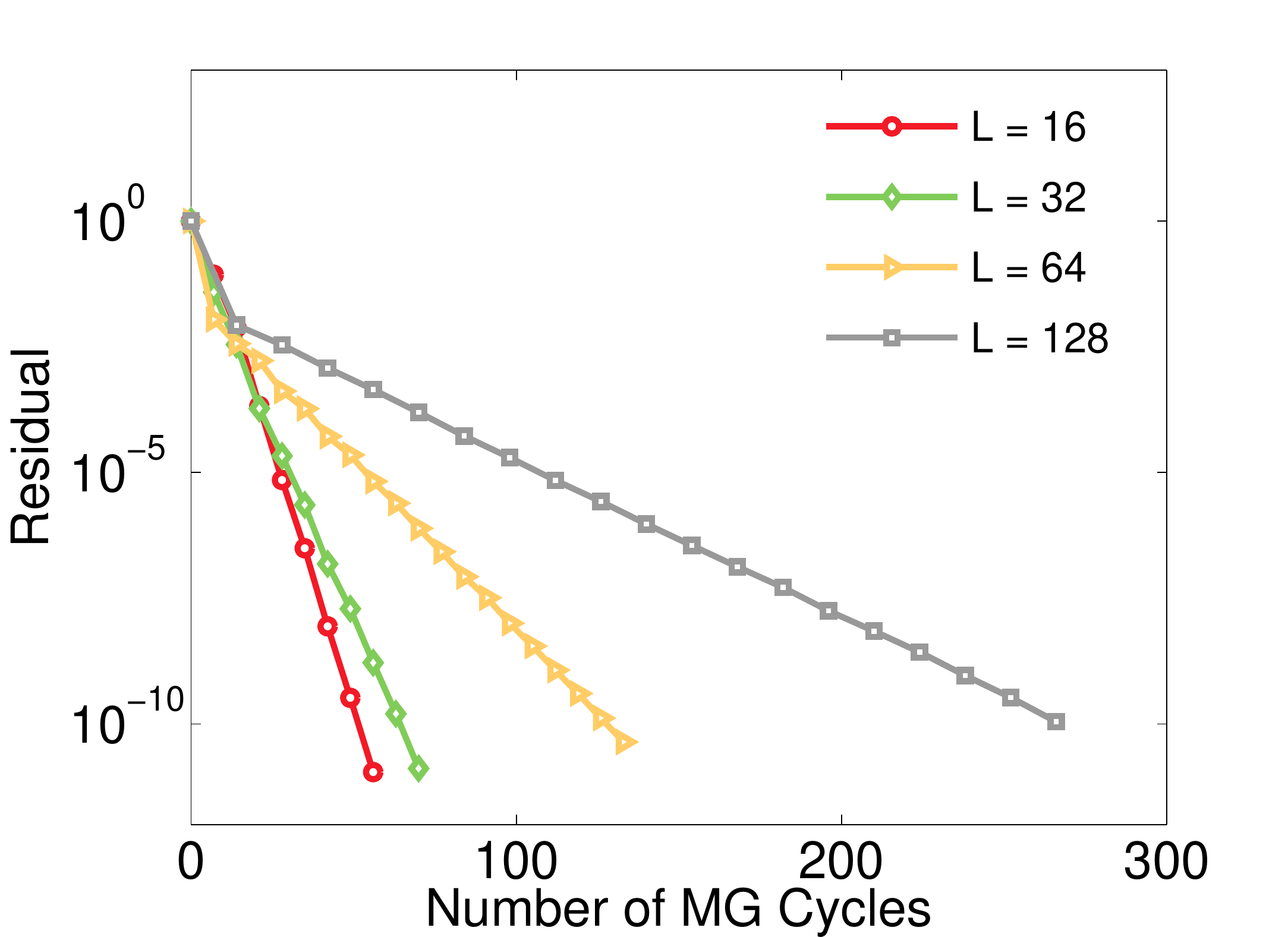}

\includegraphics[width=0.49\textwidth]{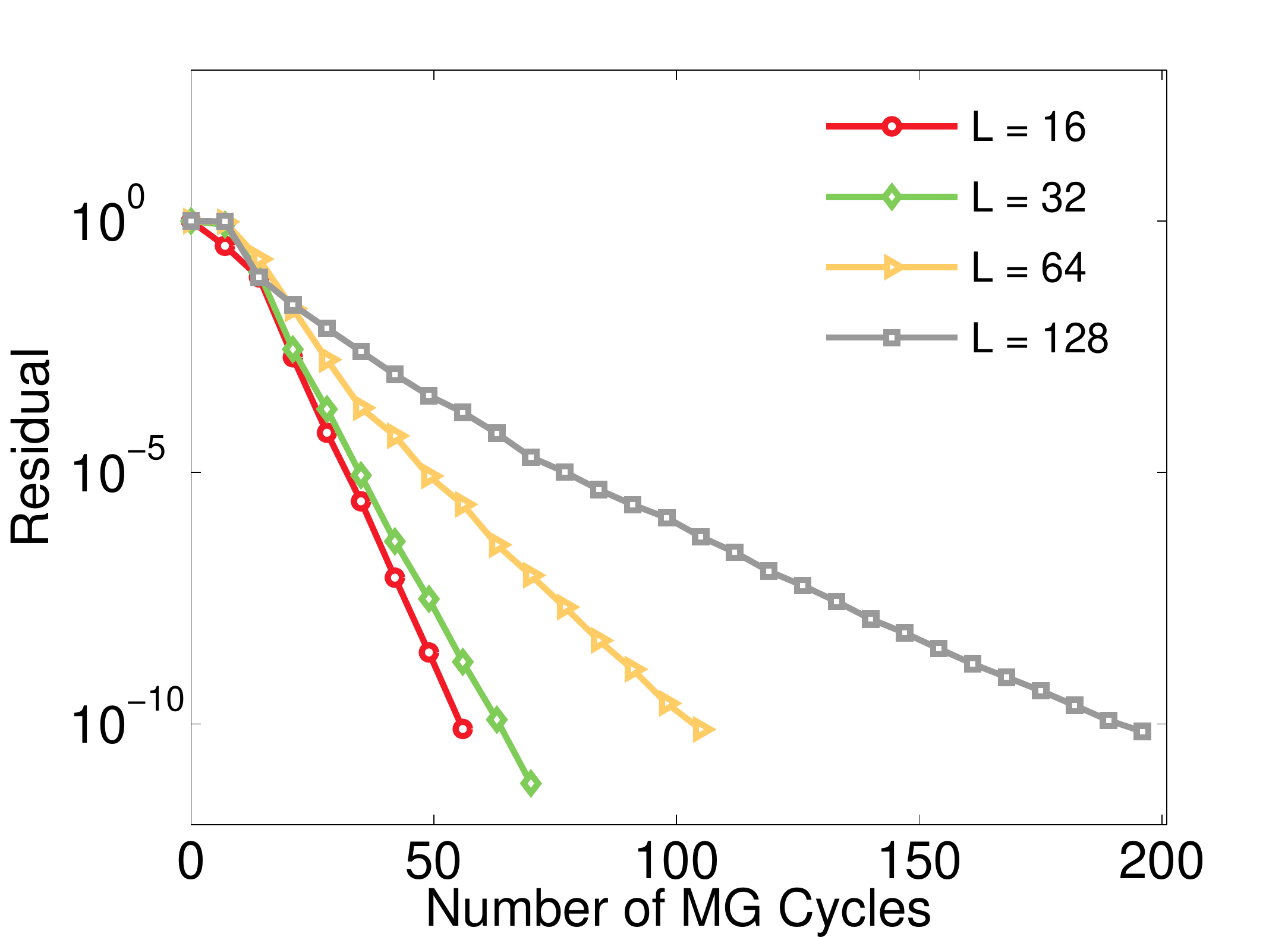}\includegraphics[width=0.49\textwidth]{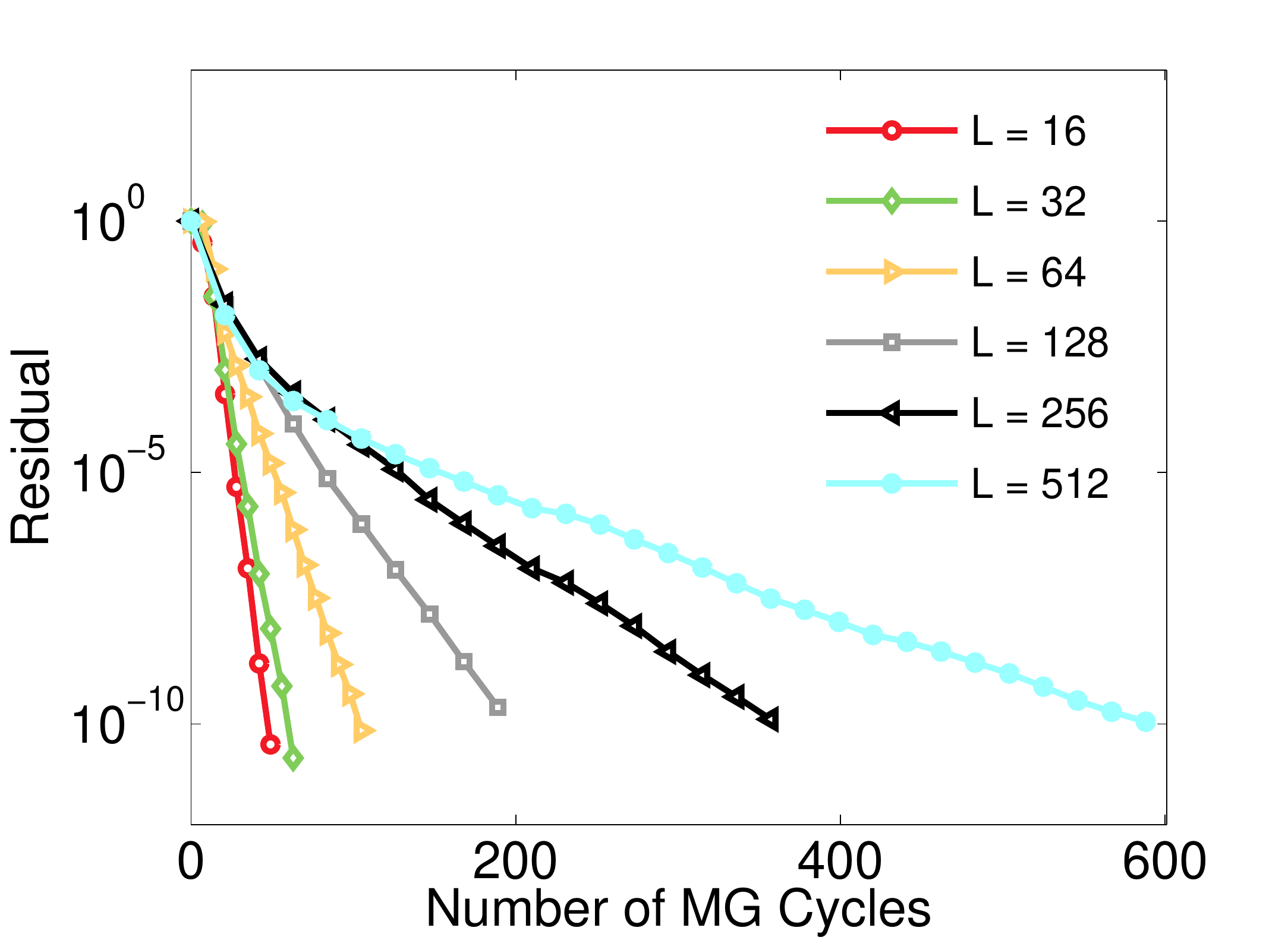}
\caption{\label{fig:InexactSchur}Performance of the approximate Schur complement
preconditioner for the saddle-point problem (\ref{eq:saddle_point}).
For the top row $\M A=L^{-2}\M I-\M L$ and for the bottom row $\M A=-\M L$.
The restart frequency is set to set to 30 in the FGMRES solver. (\emph{Top
Left}) Investigation of the optimal $m$ for $n=1$ for a three-dimensional
system of $16^{3}$ cells. (\emph{Top Right}) System size dependence
for $m=5$ and $n=1$. (\emph{Bottom left}) Performance of the special
handling of the null space of $\M A$ in three dimensions for $m=5$
and $n=1$. (\emph{Bottom right}) Same as bottom left panel but in
two dimensions.}
\end{figure}

For the approximate Schur complement preconditioner, we have found
that a large $m$ is beneficial in speeding the outer Krylov solver
convergence rate, but does not help decrease the overall computing
time (i.e., total number of multigrid cycles). This is illustrated
in the top left panel of Fig. \ref{fig:InexactSchur}. Increasing
$n$ was also not found to be beneficial (data not shown). Empirically
we observe that a reasonable choice is to use $m=5$ and $n=1$. In
the top right panel of Fig. \ref{fig:InexactSchur} we show the convergence
history of the FGMRES solver in three dimensions. Similar behavior
is seen as for the diagonal preconditioner, but the total computational
cost is somewhat increased.

In the bottom two panels of Fig. \ref{fig:InexactSchur} we show the
convergence history of the FGMRES solver we study the steady-state
case $\M A=-\M L$ in both two and three dimensions. We utilize the
approximate Schur complement modified to take into account the null
space of the discrete Laplacian, as described in Appendix \ref{sub:AppendixPrecond}.
We find similar behavior as for the case of a non-singular $\M A$,
demonstrating that the null space has been handled properly in the
preconditioner.

\begin{figure}
\includegraphics[width=0.49\textwidth]{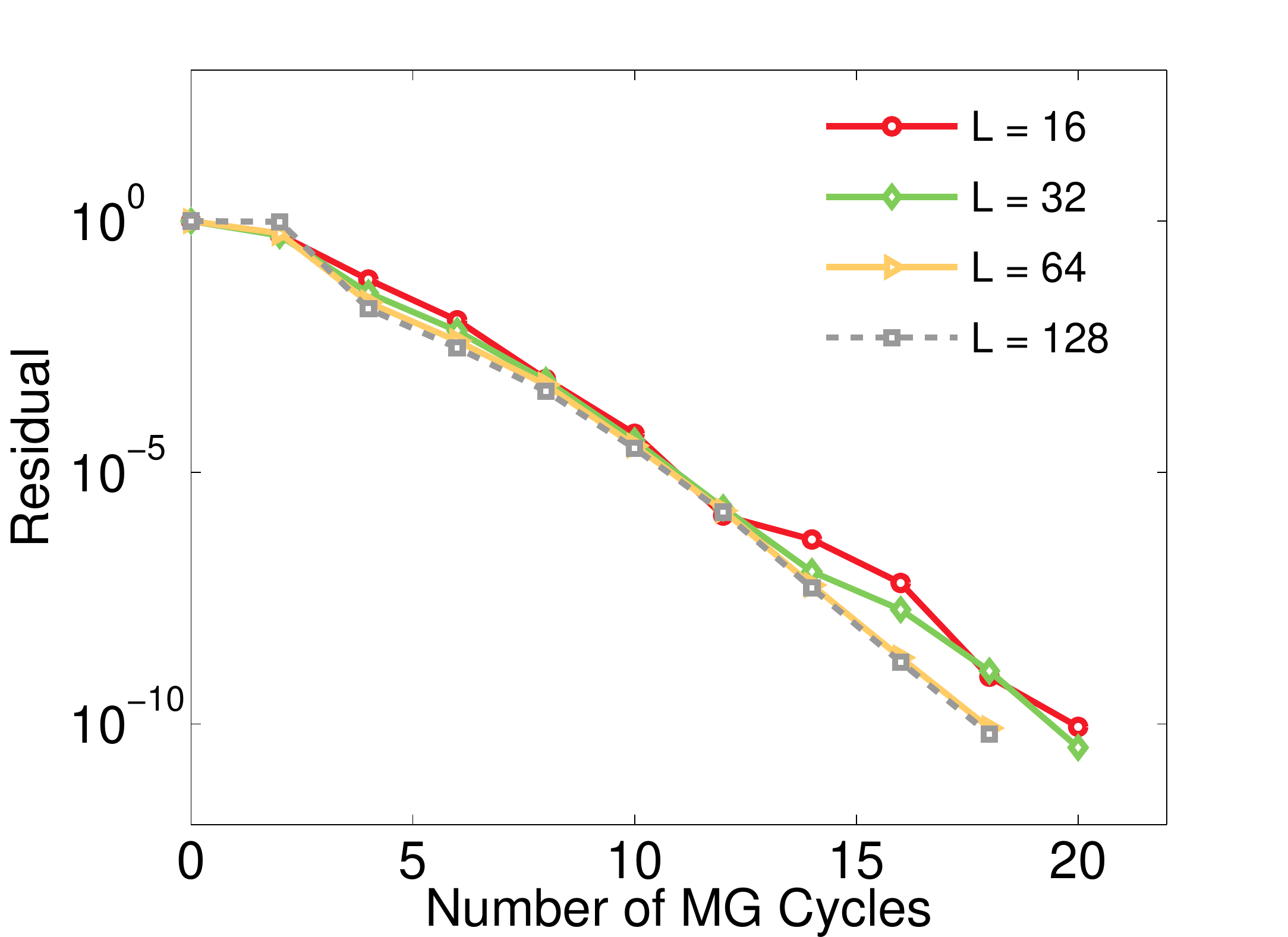}\includegraphics[width=0.49\textwidth]{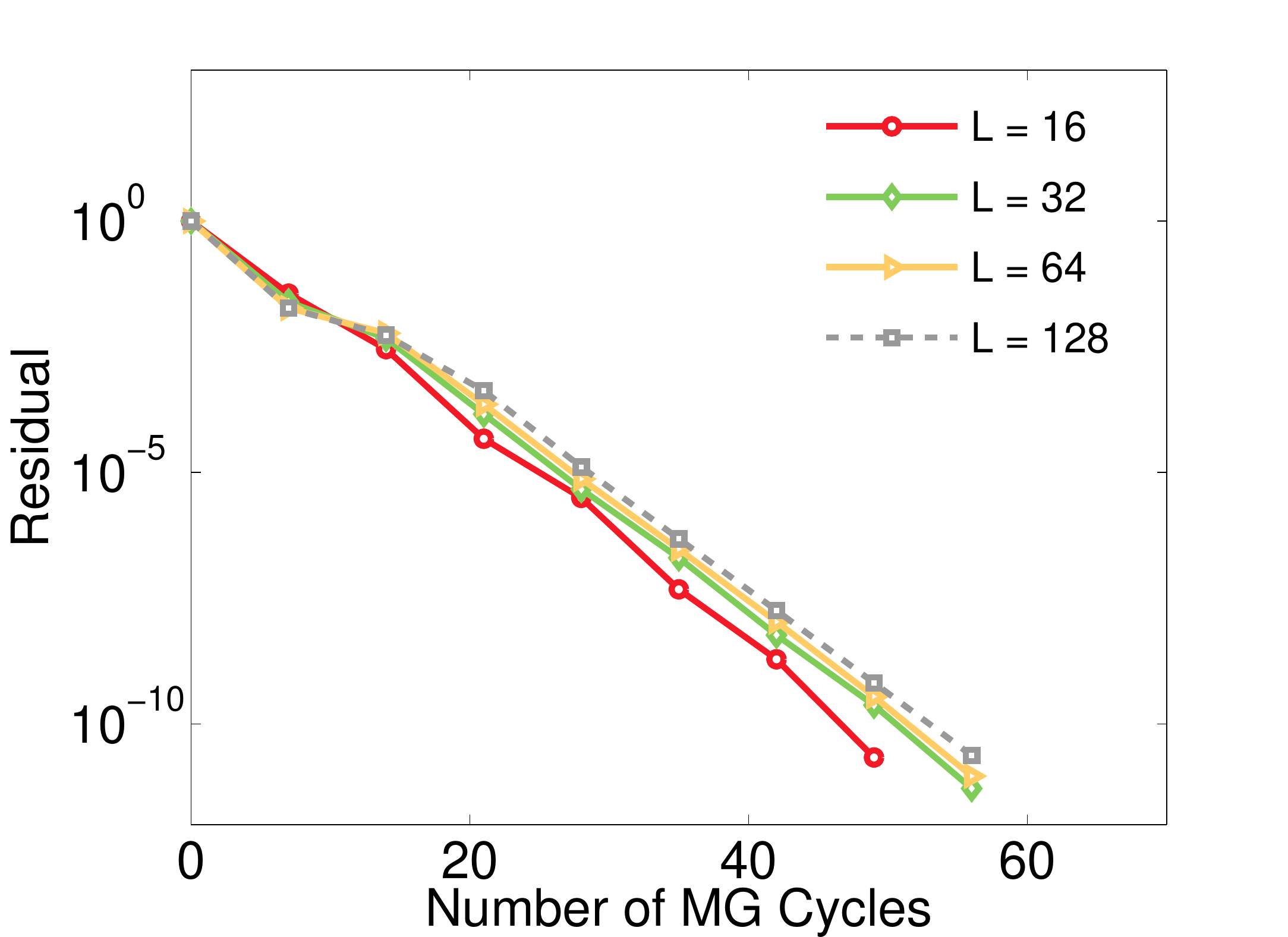}
\caption{\label{fig:FixedNScaling}Performance of the diagonal preconditioner
with $n=1$ (left panel), and the approximate Schur complement preconditioner
with $n=1$, $m=5$ (right panel) for the saddle-point problem (\ref{eq:saddle_point})
with $\M A=L^{-2}\M I-\M L$. Here the number of blobs, arranged on
a cubic lattice, is kept fixed at $64=4^{3}$ as the number of cells
$L$ along each dimension is varied.}
\end{figure}

Figs. \ref{fig:DiagonalScaling} and \ref{fig:InexactSchur} show
that neither the diagonal or the approximate Schur complement preconditioners
give a solver that is robust with respect to the system size. The
number of multigrid sweeps is seen to increase strongly as the size
of the system is increased. This is unlike multigrid solvers for the
Poisson equation, where only a mild (at most logarithmic) dependence
on system size is observed. Nevertheless, the convergence is seen
to be uniform and both preconditioners provide a viable iterative
solver, especially for small numbers of blobs. In particular, if the
number of blobs is kept small, we find that the total number of multigrid
cycles does not increase as the system size increases, as illustrated
in Fig. \ref{fig:FixedNScaling}. The approximate Schur complement
preconditioner is found to be more robust to restarts in the FGMRES
algorithm, and thus requires substantially less memory than the diagonal
preconditioner.

\section{\label{sec:Results}Results}

In this section we apply the reactive blob method to model reaction-diffusion
problems for ordered and random dispersions of reactive spheres. We
will compare the results obtained from the blob model to those obtained
by other methods \cite{ReactionDiffusion_Cubic,ReactionDiffusion_Lu,ReactionDiffusion_NDL,ReactionDiffusion_Brady,ReactionDiffusion_Torquato}
to determine the fidelity yielded by our blob model in approximating
rigid reactive spheres. A similar study was performed in the context
of fluid flow problems in Refs. \cite{IBM_Sphere,ISIBM} and it was
found that a blob can provide a surprisingly good approximation to
a rigid sphere.

\subsection{\label{sub:ReactionRadius}Reaction radius of a blob}

We expect a blob to behave similarly to a reactive sphere of radius
$a\sim h$. The constant of proportionality gives the effective reactive
radius of the blob, and depends on the discrete Laplacian $\M L$
used to discretize the Poisson equation, on the kernel function $\varphi$
used to implement the discrete local averaging and spreading operators
(see Section \ref{sub:SpatialDiscretization}), and, to a minor extent,
on the position of the blob relative to the grid of cells. Here we
calculate the effective reactive radius for the 3-point and 4-point
kernel functions \cite{IBM_PeskinReview} and the standard 7-point
discrete Laplacian in three dimensions. It is also possible to define
an effective reactive radius in two dimensions. However, in two dimensions
the behavior of a collection of reactive particles is very system
size and boundary condition-dependent, and therefore we focus here
on the more practically-relevant case of three dimensional space.

We consider a single blob in a periodic domain $\Omega$ of volume
$V_{\Omega}=L^{3}h^{3}$, which is equivalent to considering an infinite
cubic array of reactive blobs. We focus on the diffusion-limited case
$P\rightarrow0$, and study the steady state in the presence of a
uniform source of reactant of strength $s$ throughout the domain.
We solve the saddle-point problem (\ref{eq:saddle_point}) for $\xi=\zeta=1$,
$\V g=s$, and $\V f=0$ using the approximate Schur complement preconditioner
described in Appendix \ref{sub:AppendixPrecond}, with $m=5$ and
$n=1$.

\begin{figure}
\centering{}\includegraphics[width=0.49\textwidth]{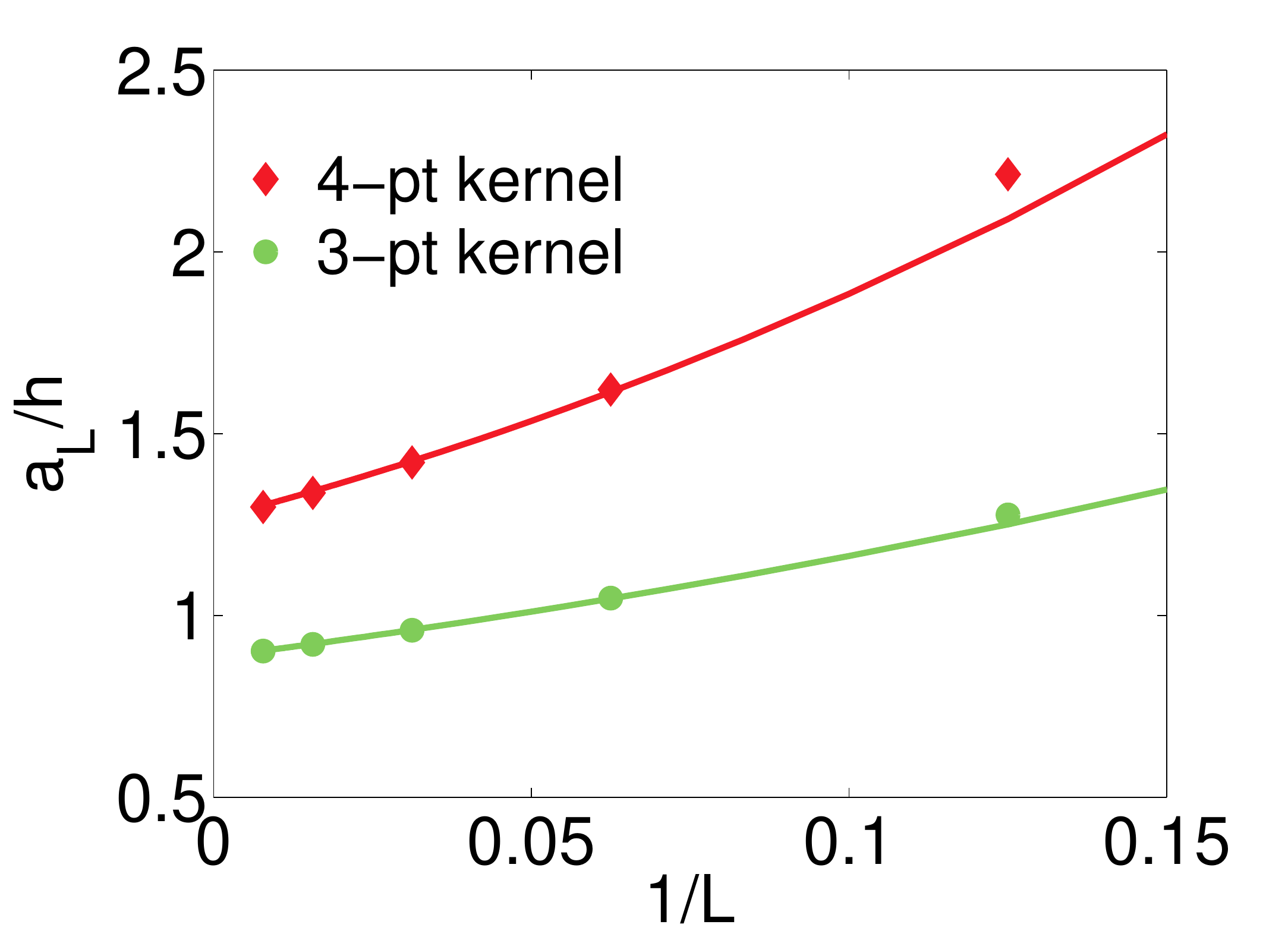}\includegraphics[width=0.49\textwidth]{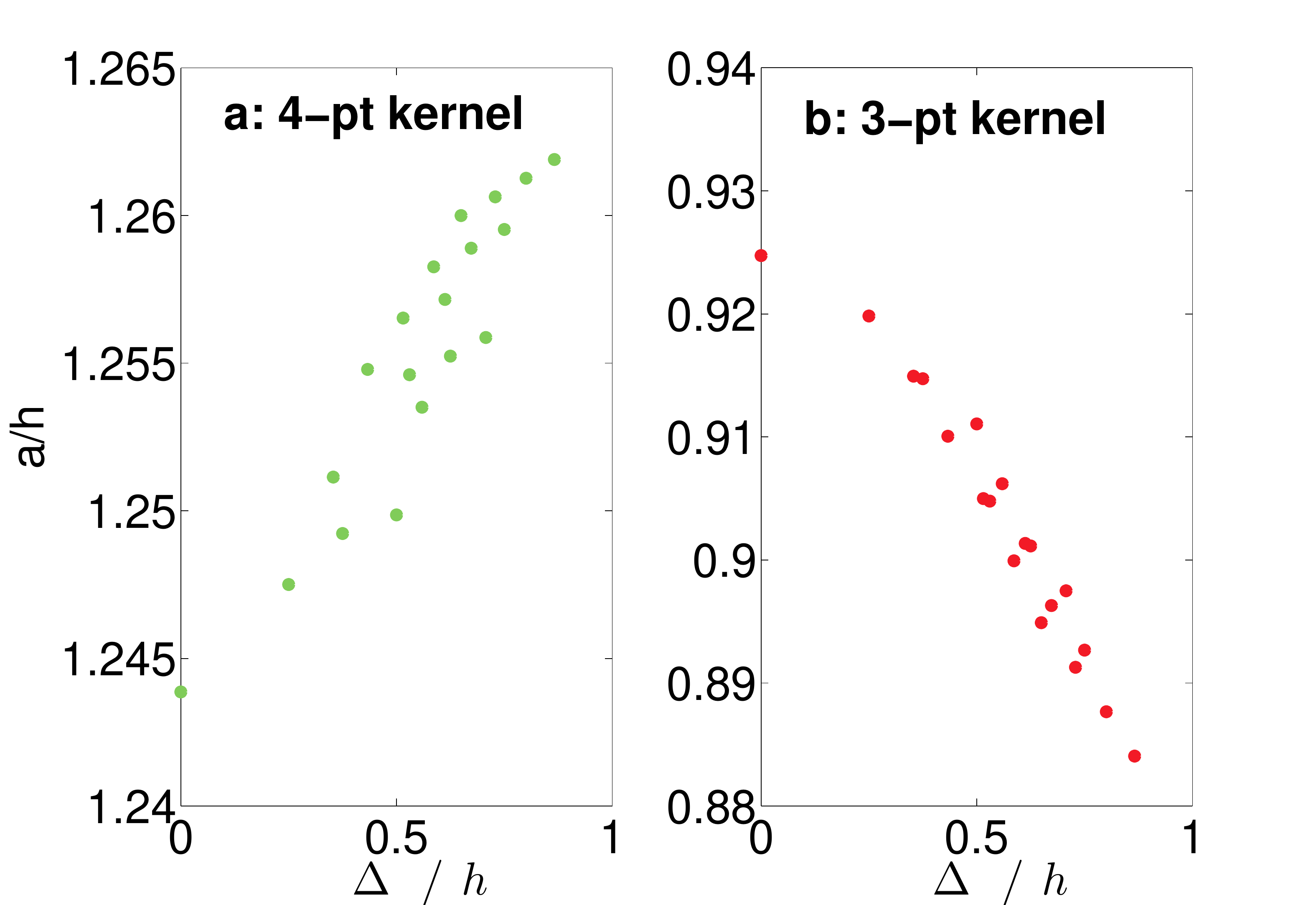}
\caption{\label{fig:ReactionRadius}Reaction radius of a blob for Peskin's
4-point and 3-point kernels. (\emph{Left}) The volume-fraction-dependent
reaction radius as a function of inverse system size, together with
fits of the form (\ref{eq:a_phi}) with $a$ as an unknown parameter.
(\emph{Right}) The translational invariance of the reactive blob radius
(extrapolated to infinite system size using (\ref{eq:a_phi})), as
the blob is displaced by a distance $\Delta$ relative to the center
of a grid cell along several different directions.}
\end{figure}

For a single reactive sphere, the normalized reaction rate (\ref{eq:beta_P_cont})
is given by 
\begin{equation}
\beta_{0}=\frac{sa^{2}}{\chi}\cdot\frac{\left(1-\varphi\right)}{3\varphi}\cdot\frac{1}{\bar{c}},\label{eq:beta_P_disc}
\end{equation}
where $\varphi=\left(4\pi a^{3}/3\right)/V_{\Omega}$ is the volume
fraction (packing density) of the cubic lattice. As $\varphi\rightarrow0$,
$\beta_{0}\rightarrow1$; therefore, by measuring the average concentration
$\bar{c}$ in the domain (recall that for a blob the domain includes
all cells, including those overlapping the blob) for very small volume
fractions, we can obtain the reactive radius $a$. If we take $s=1/V_{\Omega}$
we see that for $\varphi\ll1$, 
\begin{equation}
\beta_{0}\approx(4\pi a\chi\bar{c})^{-1}\approx1.\label{eq:beta0_low_phi}
\end{equation}
We define a volume fraction-dependent effective reactive radius $a_{\varphi}=\left[4\pi\chi\bar{c}\left(\varphi\right)\right]^{-1}$,
and measure it numerically for a blob for different system sizes.
In the left panel of Fig. \ref{fig:ReactionRadius} we show $a_{\varphi}/h$
as a function of $L^{-1}\sim\varphi^{1/3}$. Theoretical calculations
based on multipole expansions \cite{DiffusionLimited_Periodic} give
the asymptotic expansion
\begin{equation}
\beta_{0}=1+1.76\varphi^{1/3}+\left(1.76\right)^{2}\varphi^{2/3}+\text{h.o.t.},\label{eq:beta0_phi}
\end{equation}
which gives the corresponding expansion of $a_{\varphi}$, expressed
in terms of system size as 
\begin{equation}
a_{L}=a\left[1+2.84\frac{a}{Lh}+\left(2.84\frac{a}{Lh}\right)^{2}+\text{h.o.t.}\right].\label{eq:a_phi}
\end{equation}
By fitting the numerical data to this series expansion at small $\varphi$
we can obtain a numerical estimate for $a$. It is important to note
that the precise value of $a$ depends on the position of the center
of the blob relative to the underlying grid. For a blob at the node
(corner) of the grid, we numerically estimate $a\approx1.27\, h$
for the 4-point kernel and $a\approx0.885\, h$ for the 3-point kernel.

Next we explore how translationally-invariant the effective reactive
radius $a$ is. We fix the system size at $L=128$ and then move the
blob relative to the center of a grid cell in several test directions,
and estimate the infinite system size $a$ from (\ref{eq:a_phi})
as $a\approx a_{L}-2.84a_{L}^{2}/\left(Lh\right)$. In the right panel
of Fig. \ref{fig:ReactionRadius} we show $a$ as a function of the
distance $\Delta$ of the blob to the nearest cell center. We see
that $a$ varies only by about $1.5\%$ for the 4-point kernel, and
by about $5\%$ for the 3-point kernel, consistent with similar results
for flow problems \cite{DirectForcing_Balboa,ISIBM}. The 4-point
function is seen to provide much improved translational invariance,
and is therefore preferred except in cases where computational cost
is of primary concern. In future work we will study other choices
for the kernel functions.

\begin{figure}
\centering{}\includegraphics[width=0.5\paperwidth]{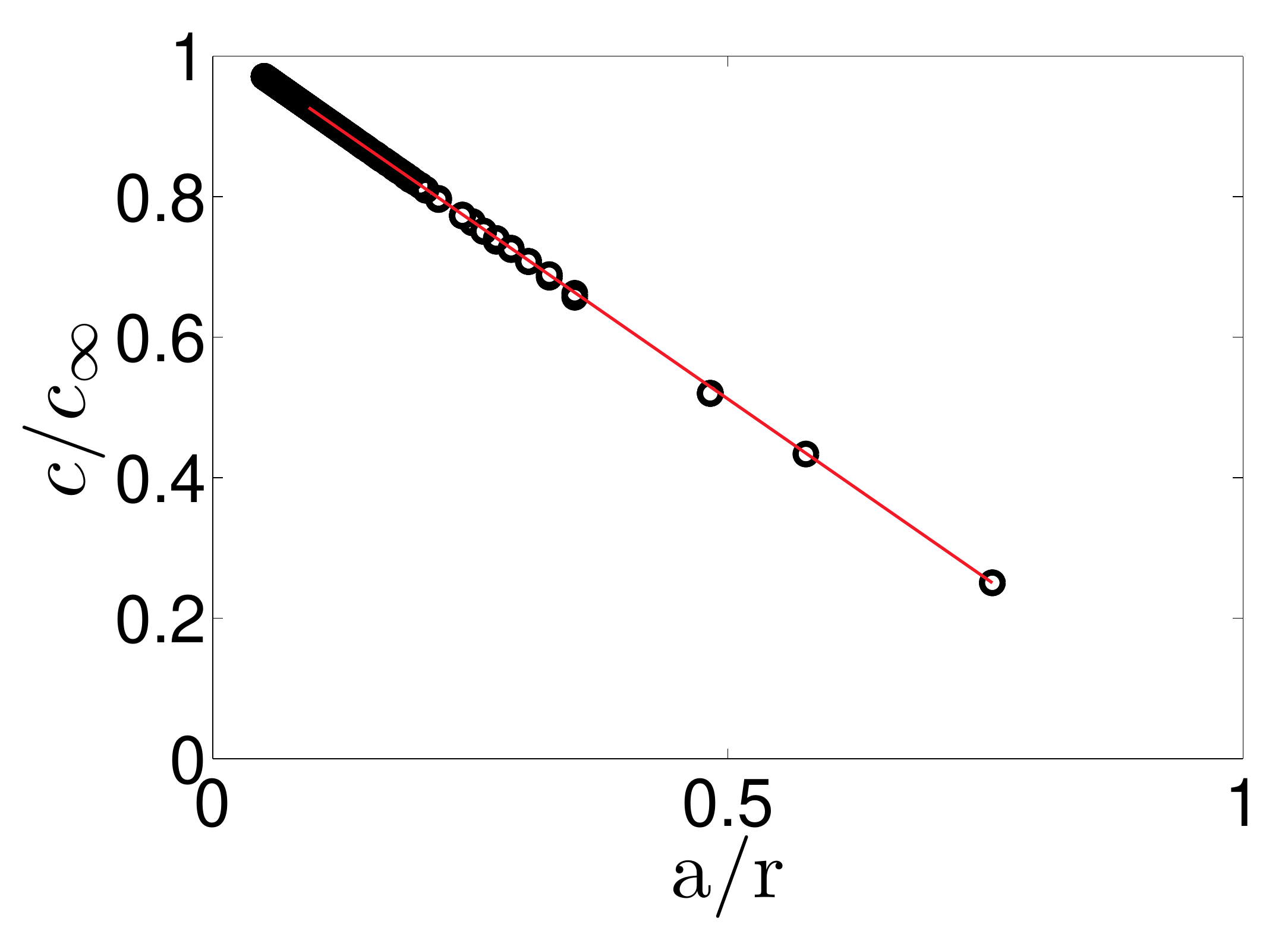}
\caption{\label{fig:FarFieldDecay}Decay of the concentration around a blob
in a cubic domain of $100^{3}$ grid cells with Dirichlet boundary
conditions, compared to the theoretical prediction (\ref{eq:c_r_approx}).
The four-point kernel $\phi_{4}$ is used so the reactive radius is
taken to be $a\approx1.255h.$}
\end{figure}

The effective reaction radius is a property of the discretized blob
(specifically, the combination of the solver for the concentration
equation and the discrete kernel function), and does not depend on
the boundary conditions. To see this, let us consider the example
of a single reactive blob in a large cubic box of size $L$ with Dirichlet
boundary conditions $c=c_{\infty}$ at the boundaries of the domain,
as first illustrated in Fig. \ref{fig:Illustration3D}. In an infinite
domain the steady-state concentration would be given by (\ref{eq:c_r_sphere}).
Here we focus on the diffusion-limited case $P=0$. We can approximate
the finite-size effects of the boundaries by solving the spherically-symmetric
problem $\grad^{2}c=0$ subject to $c(r=a)=0$ and $c(r=L/2)=c_{\infty},$
to obtain
\begin{equation}
\frac{c(r)}{c_{\infty}}\approx\frac{L}{L-2a}\left(1-\frac{a}{r}\right).\label{eq:c_r_approx}
\end{equation}
In Fig. \ref{fig:FarFieldDecay} we compare the numerical solution
for the steady-state concentration around the blob with the theoretical
result, with the value of the reactive radius $a$ taken to the average
in Fig. \ref{fig:ReactionRadius}. Excellent agreement is observed
and the $r^{-1}$ decay of the concentration around the blob is clearly
illustrated.

\subsection{Periodic Dispersions}

We now investigate how well a periodic dispersion of reactive blobs
approximates the behavior of a periodic dispersion of reactive spheres
of radius $a$. We consider a single blob in a periodic domain $\Omega$
of volume $V_{\Omega}=L^{3}h^{3}$, which is equivalent to considering
an infinite cubic array of reactive blobs.

\begin{figure}
\centering{}\includegraphics[width=0.5\textwidth]{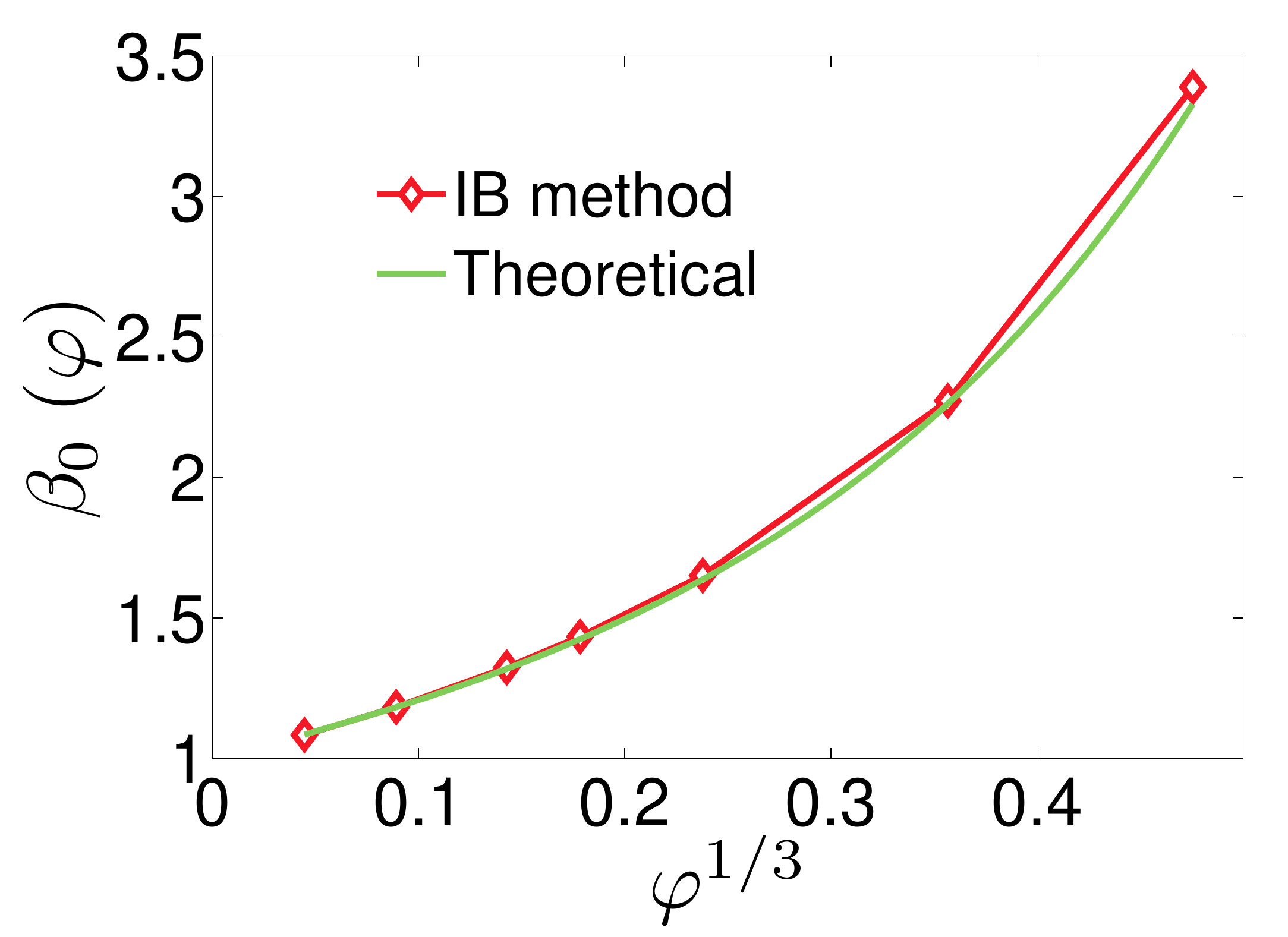}\includegraphics[width=0.5\textwidth]{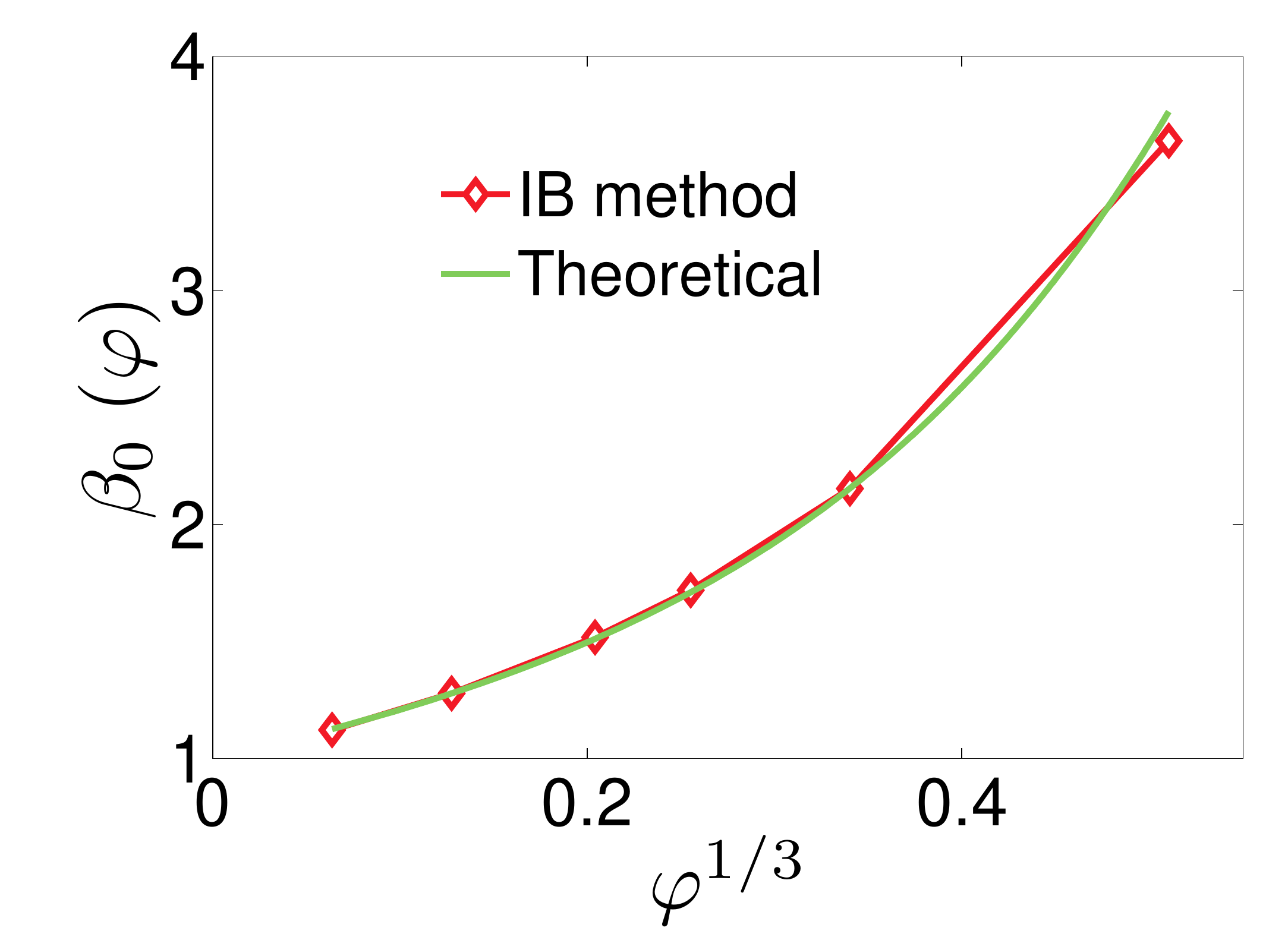}
\caption{\label{fig:PerDisp}The effective reaction rate $\beta_{0}$ for a
cubic array of blobs at several volume fractions, for the 3-point
(left panel) and the 4-point (right panel) kernels. A theoretical
prediction for a cubic array of spheres based on (\ref{eq:beta_cubic})
is shown for comparison.}
\end{figure}

In the first test we focus on the diffusion-limited case $P\rightarrow0$
and solve the steady-state (\ref{eq:saddle_point}) system in the
presence of a uniform source of reactant of strength $s$ throughout
the domain. We obtain the normalized effective reaction rate from
the average concentration in the domain using (\ref{eq:beta_P_disc}).
The numerical results are shown in Fig. \ref{fig:PerDisp}, for both
the 3-point and 4-point kernels. Based on the theoretical prediction
(\ref{eq:beta0_phi}), we fit the numerical data of Lu\emph{ }\cite{ReactionDiffusion_Cubic}
for a cubic lattice of spheres at packing fraction $\varphi=0.1$
and $\varphi=0.2$ to
\begin{equation}
\beta_{0}=1+1.76\varphi^{1/3}+\left(1.76\right)^{2}\varphi^{2/3}+b\varphi+c\varphi^{4/3},\label{eq:beta_cubic}
\end{equation}
with fitting (interpolation) parameters which we estimated to be $b\approx-0.92$
and $c\approx17.4$. We take this fit as a good approximation to the
true answer for a dispersion of spheres, and compare our numerical
data for a dispersion of blobs to this approximation in Fig. \ref{fig:PerDisp}.
An excellent match is observed over this range of $0<\varphi<0.15$.

Note that achieving higher packing densities with blobs is not possible
for a cubic arrangement without overlapping the kernels of the blobs.
This is because a cubic arrangement is a very low-density packing
of spheres. For cubic arrangements, the blob model gives unphysical
results, including negative $\beta_{0}$, at higher packing densities.
This is because inside the support of the kernel the concentration
is negative (in order to satisfy the constraint $\M{\mathcal{J}}\M c=0$,
as illustrated in Fig. \ref{fig:Illustration3D}), and averaging the
concentration over the whole domain (including the blob interiors),
can produce a negative number. It is not surprising that the minimally-resolved
blob model cannot provide a good approximation for densely-packed
spheres. How good of an approximation the blob model provides depends,
in fact, not just on the packing density but also on the type of arrangement.
In the next section we study random dispersions and find that the
blob model gives reasonably accurate results to rather high packing
densities.

\subsection{Random Dispersions}

\begin{figure}
\centering{}\includegraphics[width=0.5\textwidth]{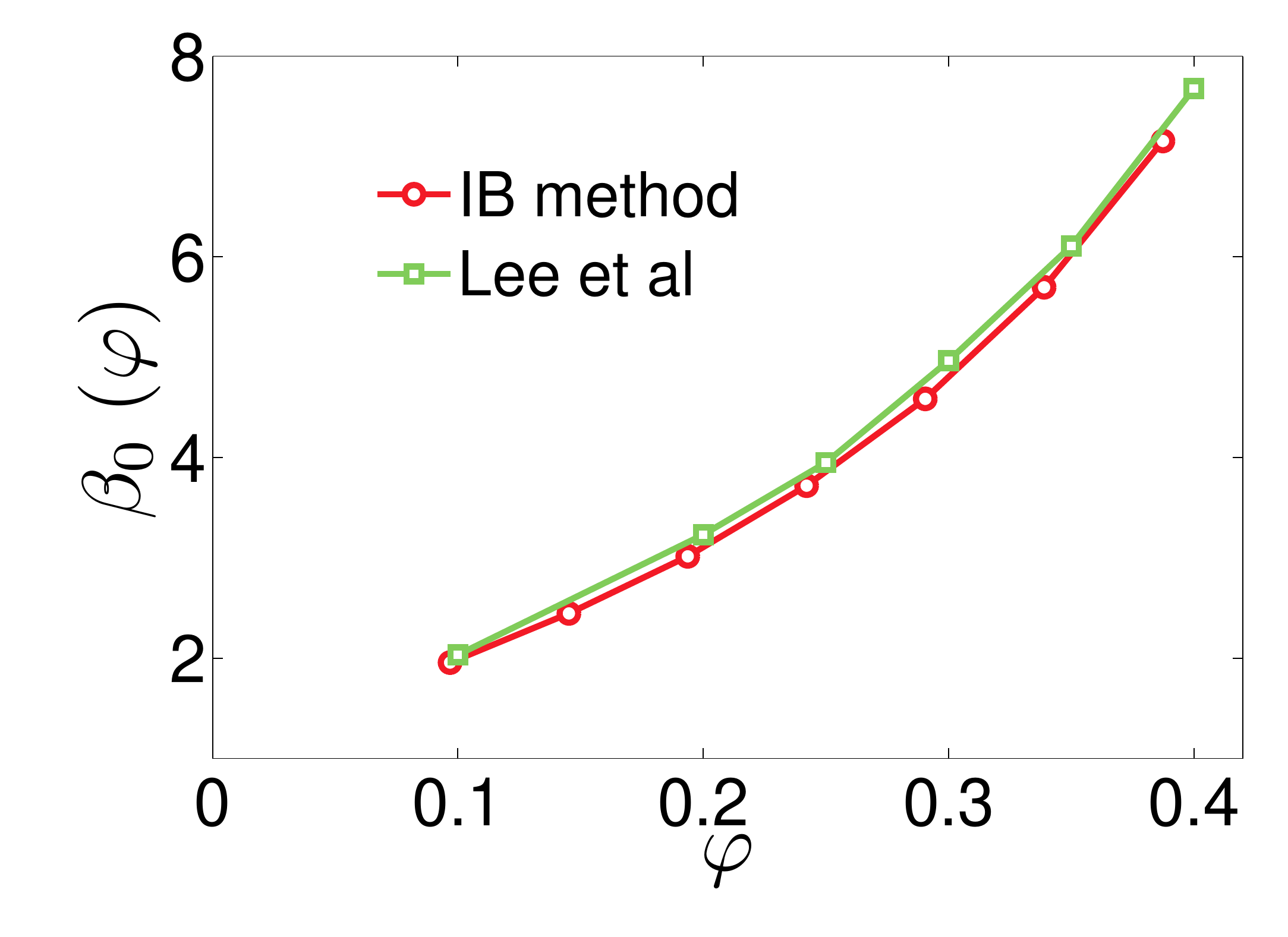}\includegraphics[width=0.5\textwidth]{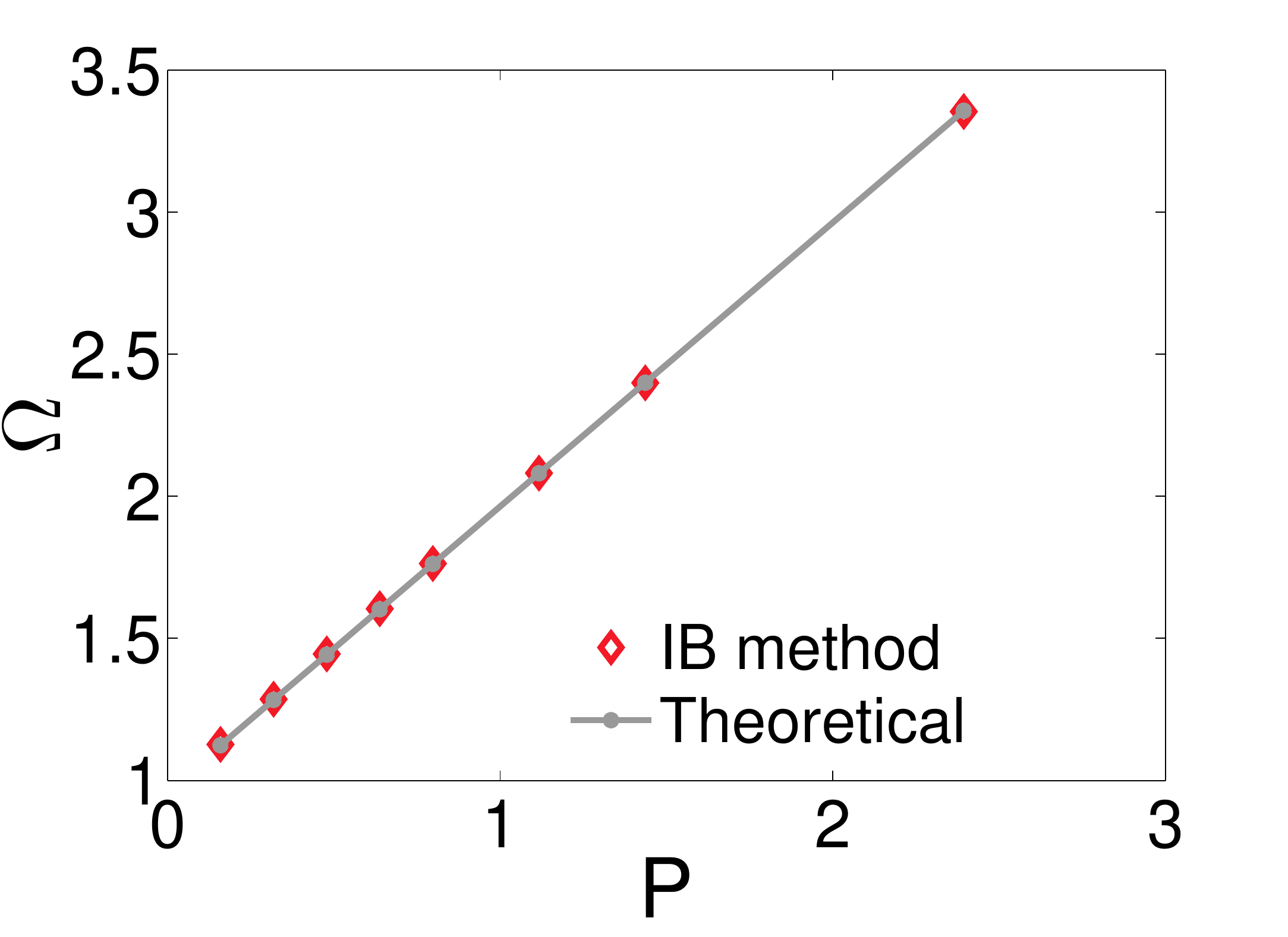}
\caption{\label{fig:randDisp}(\emph{Left}) The effective reaction rate $\beta_{0}$
for a random arrangement of non-overlapping spheres at several volume
fractions, for the 4-point kernel. The first-passage results of Lee
\emph{et al.} \cite{ReactionDiffusion_Torquato} are shown for comparison.
(\emph{Right}) \label{fig:finite_P_fit}The normalized inverse reaction
rate $\Omega$ as a function of the dimensionless number $P$ for
a single blob in a large cubic periodic box, compared with the theoretical
prediction $\Omega=\beta_{0}^{-1}+(1-\varphi)P$.}
\end{figure}

We now turn to random dispersions of blobs and compare to the first-passage
results of Lee \emph{et al.} \cite{ReactionDiffusion_Torquato} for
random dispersions of non-overlapping spheres. We generate configurations
of non-overlapping hard spheres at several volume fractions using
the packing algorithm and code described in Refs. \cite{MRJ_HS_4D,Jamming_g2}.
The computational grid is kept at $L^{3}=128^{3}$ cells, and the
number of particles is chosen to give (approximately) a desired packing
fraction in the range $0.1\leq\varphi\leq0.4$ (the densest dispersion
had on the order of $10^{5}$ particles). Blobs are placed at the
center of each hard sphere, and the steady-state diffusion-limited
equation (\ref{eq:saddle_point}) is solved numerically. In the left
panel of Fig. \ref{fig:randDisp} we compare the data for blobs to
the data for rigid spheres \cite{ReactionDiffusion_Torquato}. An
excellent agreement is observed even for packing fractions as high
as $\varphi=0.4$, surprisingly close to the jamming density $\varphi\approx0.64$.
This unexpected accuracy is perhaps owing to cancellation of errors
due to the randomness. Note that at higher packing densities a large
fraction of the blob kernels are overlapping even though the equivalent
rigid spheres would not be.

\subsection{\label{sub:FiniteP}Finite reaction rate}

The results of the previous section were obtained for diffusion-limited
case $P\rightarrow0$. Here we consider the same setup of a single
blob in a periodic domain, but consider the steady state for a finite
reaction rate $k$. Using (\ref{eq:kappa_k}), the dimensionless number
$P=\chi/\left(ka\right)$ for a sphere of radius $a$, becomes $P=4\pi\chi a/\kappa$
for a blob, where $a$ is the reactive radius of the blob. 

Lu \cite{ReactionDiffusion_finiteP} has shown that at small volume
fractions the effective reaction rate for a finite $P$ can be related
to the reaction rate for $P=0$ via a simple approximation,
\begin{equation}
\beta_{P}=\frac{1+P}{\Omega}\approx\frac{(1+P)}{\beta_{0}^{-1}+(1-\varphi)P},\label{eq:beta_P_theory}
\end{equation}
where $\varphi$ is the volume fraction and $\beta_{0}$ is the normalized
reaction rate for the same configuration of the spheres in the diffusion-limited
case. Here 
\begin{equation}
\Omega=\left[\frac{sa^{2}}{\chi}\cdot\frac{\left(1-\varphi\right)}{3\varphi}\cdot\frac{1}{\bar{c}}\right]^{-1}\approx\beta_{0}^{-1}+(1-\varphi)P\label{eq:Omega}
\end{equation}
can directly be measured from the steady state average concentration
in the domain. Equations (\ref{eq:beta_P_theory},\ref{eq:beta0_phi})
give us a theoretical prediction for $\beta_{P}$ for a cubic array
of spheres at low packing density.

We numerically solve for the steady-state concentration around a single
reactive blob at several finite reaction rates, using the 4-point
kernel. We take a box of size $L=100$ grid cells, which makes the
finite-size effects rather small. The linear relationship (\ref{eq:Omega})
between $\Omega$ and $P$ is shown in the right panel of Fig. \ref{fig:finite_P_fit}
and seen to be in excellent agreement with the numerical data over
the whole range of $P$ values, confirming that the blob model captures
the effect of finite reaction rate accurately.

\begin{figure}
\centering{}\includegraphics[width=0.5\textwidth]{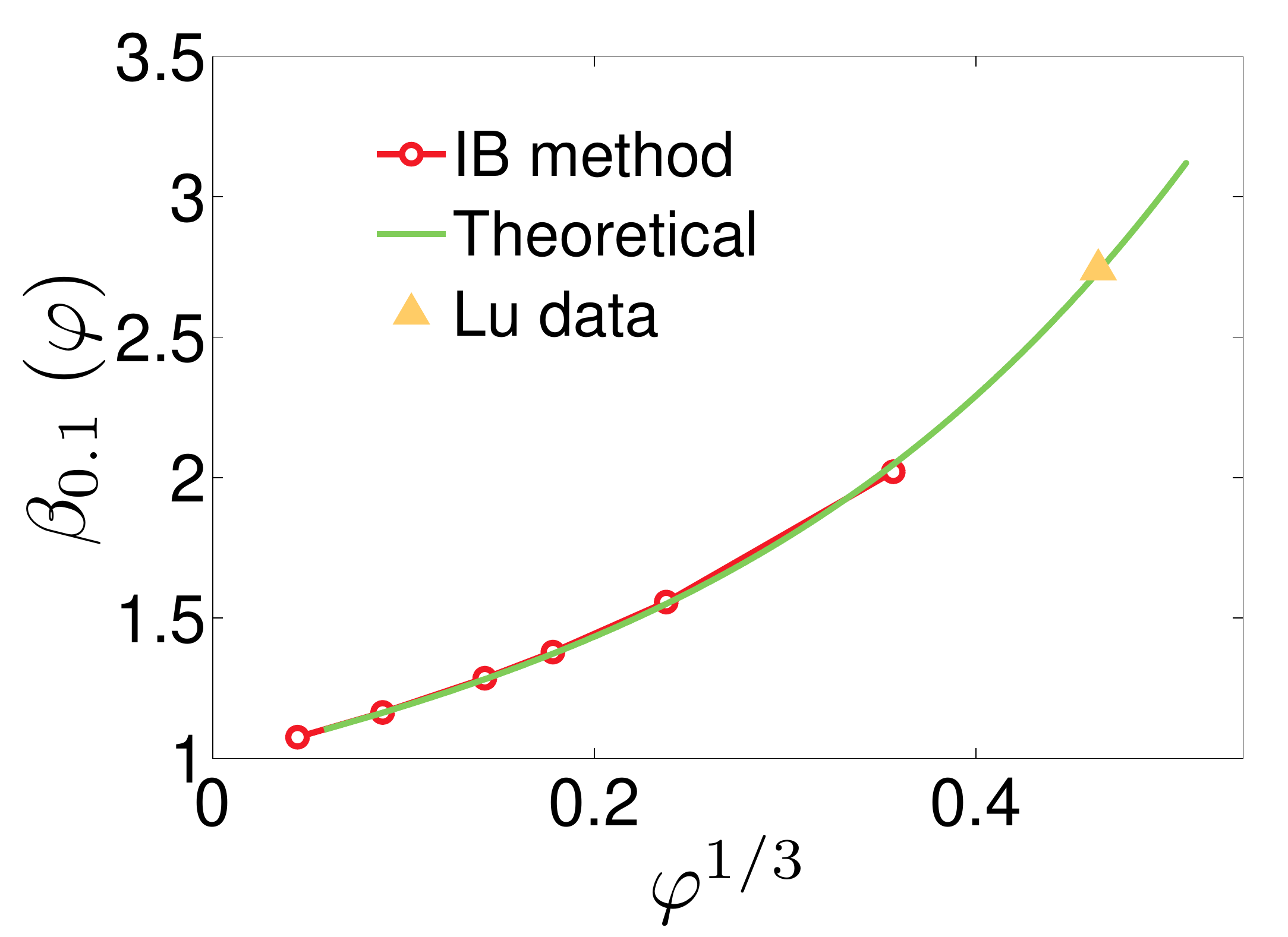}\includegraphics[width=0.5\textwidth]{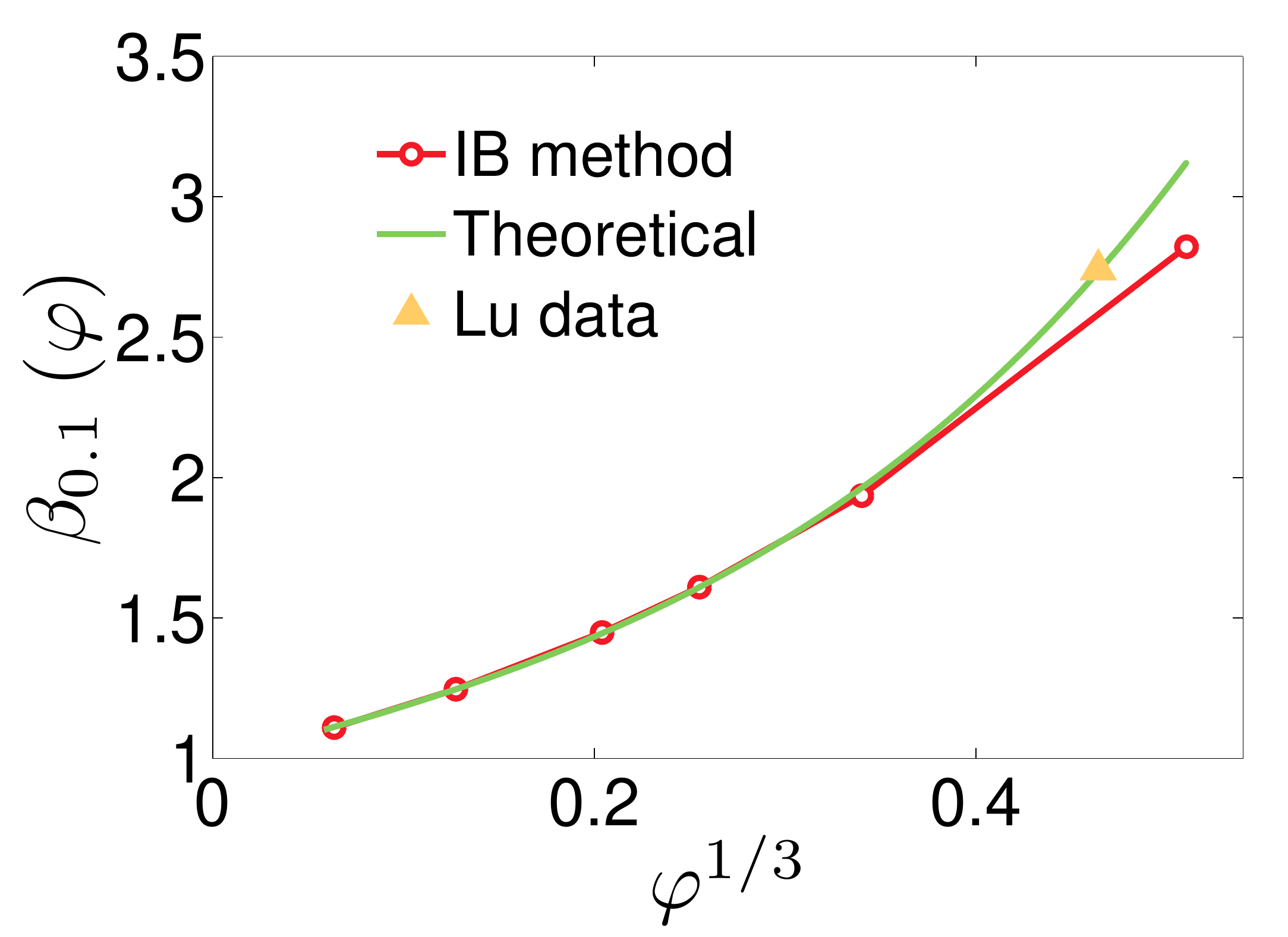}
\caption{\label{fig:finite_P_01}The effective reaction rate $\beta_{0.1}$
for a cubic array of spheres at several volume fractions, for the
3-point (left panel) and the 4-point (right panel) kernels. A theoretical
prediction based on (\ref{eq:beta_P_theory},\ref{eq:beta_cubic})
is shown for comparison. }
\end{figure}

In Fig. \ref{fig:finite_P_01} we examine the particular case of finite
reaction rate, $P=0.1$ for a cubic lattice of reactive blobs. We
compare the blob results to the theoretical prediction (\ref{eq:beta_P_theory},\ref{eq:beta_cubic})
for several packing fractions. A very good agreement is observed,
although with some clear deviations from the theory and the numerical
data of Lu \cite{ReactionDiffusion_Cubic} at the larger packing fractions.
These results confirm that the blob model provides an effective approximation
at low and moderate packing densities.

\section{Conclusions}

We developed a diffuse particle approach to modeling reaction-diffusion
processes in particle dispersions that is based on the immersed boundary
method \cite{IBM_PeskinReview} and can be used from the diffusion-limited
to the reaction-limited setting. Building on positive experience with
recently-developed methods for low-density colloidal suspensions \cite{DirectForcing_Balboa,ISIBM},
we represented each reactive particle with a diffuse object termed
a blob. In this minimally-resolved approach, each three-dimensional
blob interacts with less than $\sim50$ grid points ($3^{3}=27$ for
the three-point kernel and $4^{3}=64$ for the four-point kernel),
unlike more traditional boundary-integral, traditional immersed-boundary
methods, or other discretizations in which the particle surface is
explicitly resolved.

A reactive blob, just like a reactive sphere, is characterized only
by a reactive radius. This reactive radius is proportional to the
width of the kernel function used in the immersed boundary handling
of the interaction between the blob and the computational grid used
to solve Poisson's equation for the reactant concentration. For the
standard kernel functions employed here the reactive radius is proportional
to and, in fact, close to the grid spacing. In the future, we will
explore the construction of kernels with variable width, as necessary
to represent polydisperse dispersions of particles. Adaptive mesh
refinement can be used in cases of very large polydispersity, or in
cases when the density of blobs is strongly heterogeneous.

By comparing to accurate results in the literature, obtained using
computationally expensive Monte Carlo techniques, we demonstrated
numerically that the blob model can provide an accurate representation
at low to moderate packing densities of the reactive spheres, at a
cost similar to solving several Poisson equations in the same domain.
The blob model was shown to be effective at both finite reaction rates
and the diffusion-limited case of infinite reaction rate.

At the level of the formulation, the blob method is very similar to
the method of multipole expansions \cite{ReactionDiffusion_Brady},
truncated at the leading-order (monopole) level and regularized to
avoid singularities. An essential difference is that our method does
not require analytically-computed Green's functions, but rather effectively
computes regularized discrete Green's functions ``on the fly'' by
using a standard grid-based discretization of the Poisson equation.
This allows for great flexibility in implementing different boundary
conditions, coupling to fluid flow (including immersed-boundary methods
for flow problems), and the inclusion of other effects such as temporal
evolution, nonlinearities, and other terms not easily handled in a
Green's function-based approach.

More complicated particle shapes can be built out of a collection
of reactive spheres, as done, for example, in Refs. \cite{IBM_Sphere,StokesianDynamics_Rigid}
for flow problems around rigid bodies. In fact, a more accurate representation
of a spherical particle can be achieved by representing it with a
collection of blobs (regularized monopoles). An alternative approach
is to include dipole terms in the blob description by adding additional
degrees of freedom to each blob \cite{ReactionDiffusion_Brady}. Such
an approach has been proposed by Maxey and collaborators in the context
of fluid flow problems \cite{ForceCoupling_Stokes}, and can be extended
to reaction-diffusion problems in the diffusion-limited regime. At
the dipole level, in addition to the unknown sink strength $\lambda$,
one would associates with each particle an unknown dipole strength
(a vector) $\M{\Lambda}$, which corresponds to the continuum dipole
moment \cite{ReactionDiffusion_Lu}
\[
\M{\Lambda}=\chi\int_{\partial\mathcal{B}}\V n\cdot\left(\V r\otimes\grad c-c\M I\right)dS.
\]
In the blob formulation $\M{\Lambda}$ would be a Lagrange multiplier
associated with a ``rigidity'' constraint that imposes that the
concentration be constant (alternatively, that the concentration gradient
be zero) in the interior of the blob \cite{ForceCoupling_Stokes}.
Solving the extended linear system for $\V c$, $\lambda$ and $\V{\Lambda}$
would require the development of novel linear-algebra techniques.

At present, linear algebra is a remaining bottleneck in scaling blob
calculations to very large number of particles in the diffusion-limited
regime. In this work we proposed and studied simple preconditioners
based on using standard geometric multigrid techniques for Poisson
problems as a preconditioner to a Krylov iterative solver for the
saddle-point problem (\ref{eq:saddle_point}). This was found to work
well for small numbers of blobs, but the conditioning number increases
as the number of blobs increases, a signature of the long-ranged interactions
between the blobs. An alternative approach, which we hope to explore
in the future, is to use recently-developed specialized geometric
multigrid solvers for implicit immersed-boundary methods \cite{IBMultigrid_Guy,GeometricMultigrid_Boyce}.
In the diffusion-limited regime, it may be possible to use these multigrid
solvers with a suitably chosen finite value of $\kappa$ as a preconditioner
for the saddle-point system. It remains to be determined how effective
such preconditioners are in solving linear systems arising in the
commonly-occurring case of particles interacting through a continuum
via long-ranged power-law interactions.
\begin{acknowledgments}
A. Donev was supported in part by the Office of Science of the U.S.
Department of Energy through Early Career award DE-SC0008271 and by
the Air Force Office of Scientific Research under grant number FA9550-12-1-0356.
B. Griffith acknowledges research support from the National Science
Foundation under awards OCI 1047734 and DMS 1016554. A. Bhalla and
N. Patankar acknowledge research support from the National Science
Foundation (NSF awards CBET-0828749, CBET-1066575, and CMMI-0941674
awarded to N. Patankar). Computational resources were provided by
Northwestern University High Performance Computing System \textendash{}
Quest.
\end{acknowledgments}
\begin{appendix}

\section{\label{sub:AppendixPrecond}Preconditioning the Saddle-Point Solver}

In this Appendix we give a detailed description of the algorithm used
to precondition the saddle-point system (\ref{eq:saddle_point}).
The preconditioner can be thought of as an approximate solver for
that linear system.

For a periodic system at steady state, $\M A=-\chi\M L$ is not invertible
because all vectors in the range of $\M A$ have zero average value,
$\av{\M A\V c}=0$ for any $\V c$, where $\av{\cdot}$ denotes a
spatial average over the domain, $\av{\V c}=N_{c}^{-1}\sum_{k=1}^{N_{c}}c_{i}$,
and $N_{c}$ is the number of cells in the grid. In this case we need
to modify the Schur complement approximation. From the fact that $\M A\V c$
has zero mean, we know that
\[
\zeta\av{\M{\mathcal{S}}\V{\lambda}}=\av{\V g},
\]
Spreading preserves the sum of values because of the normalization
condition on the kernel function, 
\[
\av{\M{\mathcal{S}}\V{\lambda}}=N_{c}^{-1}\D V_{f}^{-1}\sum_{i=1}^{N}\lambda_{i}=N\, V_{\Omega}^{-1}\av{\V{\lambda}},
\]
where $\av{\V{\lambda}}$ denotes an average over all blobs, $N$
is the number of blobs, $\D V_{f}$ is the volume of a grid cell,
and $V_{\Omega}$ is the volume of the domain. This means that an
additional condition for solvability is
\[
\av{\V{\lambda}}=V_{\Omega}\zeta^{-1}N^{-1}\av{\V g},
\]
which is to be added to the saddle-point problem (\ref{eq:saddle_point})
as an additional constraint. This augmented saddle-point problem has
a unique solution and can be solved using the FGMRES algorithm.

Let us set $\alpha=1$ when $\M A$ has a null-space, and $\alpha=0$
when $\M A$ is invertible. When $\alpha=1$, let us define the restricted
inverse $\M A^{-1}$ to only act on vectors of mean value zero, and
to return a vector of mean zero. The approximate Schur complement
preconditioner and diagonal preconditioners can be implemented as
follows: 
\begin{enumerate}
\item Compute a concentration estimate 
\[
\tilde{\V c}^{\star}=\widetilde{\M A}_{n}^{-1}\tilde{\V g}=\widetilde{\M A}_{n}^{-1}\left(\V g-\alpha V_{\Omega}N^{-1}\M{\mathcal{S}}\av{\V g}\right),
\]
where the notation $\M{\mathcal{S}}q$ denotes $\M{\mathcal{S}}$
applied to a vector of length $N$ with all entries equal to $q$.
Here $\widetilde{\M A}_{n}^{-1}$ is an approximation of $\M A^{-1}$
obtained by using $n$ cycles of geometric multigrid, starting from
a zero initial guess.
\item Solve the Schur complement system
\begin{equation}
\left(\M{\mathcal{J}}\widetilde{\M A}^{-1}\M{\mathcal{S}}\right)\tilde{\V{\lambda}}=\left(\zeta\xi\right)^{-1}\tilde{\V h}=\left(\zeta\xi\right)^{-1}\left(\xi\M{\mathcal{J}}\tilde{\V c}^{\star}-\V f\right),\label{eq:approx_Schur}
\end{equation}
approximately. If $\alpha=1$, ensure that $\av{\tilde{\V{\lambda}}}=0$
in the end by subtracting the mean.

\begin{enumerate}
\item For the diagonal preconditioner, the Lagrange multipliers would be
estimated as
\[
\tilde{\V{\lambda}}\approx\chi\left(\zeta\xi\gamma\right)^{-1}\tilde{\V h}.
\]

\item For the approximate Schur complement preconditioner, use $m$ iterations
of a Krylov method to solve (\ref{eq:approx_Schur}), approximating
$\widetilde{\M A}^{-1}=\widetilde{\M A}_{1}^{-1}$ by using a single
cycle of geometric multigrid.
\end{enumerate}
\item Correct the concentration
\[
\tilde{\V c}=\widetilde{\M A}_{n}^{-1}\left(\V g-\zeta\M{\mathcal{S}}\tilde{\V{\lambda}}\right),
\]
using the previous estimate $\tilde{\V c}^{\star}$ as an initial
guess for the multigrid solver.
\item If $\alpha=1$, determine the average concentration $\bar{c}=\av c$
from the condition $\av{\xi\M{\mathcal{J}}\tilde{\V c}}=\av{\V f-\xi\bar{c}}$,
\[
\bar{c}=\av{\xi^{-1}\V f-\M{\mathcal{J}}\tilde{\V c}}.
\]

\item Return the approximation to the solution of the saddle-point problem,
\begin{align*}
\V{\lambda} & =\tilde{\V{\lambda}}+\alpha V_{\Omega}\zeta^{-1}N^{-1}\av{\V g}\\
\V c & =\tilde{\V c}+\alpha\bar{c}.
\end{align*}

\end{enumerate}
\end{appendix}


\end{document}